\documentclass{nyj}

\usepackage[mathscr]{eucal}

\title[Bounded generation of\/ $\SL(n,A)$]
 {Bounded generation of $\SL(n,A)$
 \\ (after D.~Carter, G.~Keller, and E.~Paige)}

\author{Dave Witte Morris}
\address{Department of Mathematics and Computer Science \\
 University of Lethbridge \\
 Lethbridge, Alberta T1K~3M4, Canada}
\email{Dave.Morris@uleth.ca, 
http://people.uleth.ca/$\sim$dave.morris/}



\newcommand{\ints}{\mathcal{O}}
\newcommand{\II}{\mathcal{Q}}
\newcommand{\Iq}{\mathfrak{q}}

\newcommand{\gen}[1]{\left\langle #1\right\rangle}
\newcommand{\gennum}[2]{\left\langle#1\right\rangle_{\mkern-4mu#2}}

\newcommand{\Id}[1]{\mathbb{I}_{#1 \times #1}}
\newcommand{\lang}{\mathcal{L}}
\newcommand{\theory}{\mathcal{T}}
\newcommand{\genset}{\mathcal{X}}
\newcommand{\lowset}{\mathcal{M}}
\newcommand{\conje}[1]{#1^\#}
\newcommand{\localize}[2]{\mathord{{#1}{#2}^{-1}}}
\newcommand{\BS}{\localize{B}{S}}
\newcommand{\SRquot}{\mathord{\SR_{1\frac{1}{2}}}}
\renewcommand{\mod}{\mathrel{\rm mod}}

\DeclareMathOperator{\SL}{SL}
\DeclareMathOperator{\GL}{GL}

\DeclareMathOperator{\LU}{LU}
\DeclareMathOperator{\EXP}{\mathsf{Exp}}
\DeclareMathOperator{\GEN}{\mathsf{Gen}}
\DeclareMathOperator{\UNIT}{\mathsf{Unit}}
\DeclareMathOperator{\CONJ}{\mathsf{Conj}}
\DeclareMathOperator{\SR}{\mathsf{SR}}
\DeclareMathOperator{\Elem}{E}
\DeclareMathOperator{\Vas}{SSL}
\newcommand{\Enorm}{\Elem^{\normal}}
\newcommand{\LUnorm}{\LU^{\normal}}
\DeclareMathOperator{\Norm}{N}
\DeclareMathOperator{\lcm}{lcm}

\newcommand{\rational}{\mathbb{Q}}
\newcommand{\integer}{\mathbb{Z}}
\renewcommand{\natural}{\mathbb{N}}
\newcommand{\complex}{\mathbb{C}}
\newcommand{\normal}{\triangleleft}
\newcommand{\iso}{\cong}

\newcommand{\U}[1]{U \mkern -4mu \left( #1 \right)}

\newcommand{\bigset}[2]{\left\{\, #1 
 \mathrel{\left| \vphantom {\left\{ #1 \mid #2 \right\} }
 \right.} #2 \,\right\} }

\newcommand{\Men}[3]
 {\left[\genfrac{}{}{0pt}{}{\textstyle \vphantom{b}#1}{\textstyle
#2}\right]_{#3}} \newcommand{\mat}[4]{\begin{bmatrix} #1 & #2 \\ #3 &
#4\end{bmatrix}}

\newcommand{\UMen}[3]
 {\setbox0\hbox{$\genfrac{}{}{0pt}{}
 {\vphantom{b}\textstyle #1}{\textstyle#2}$}
 \left[\vphantom{\copy0}\right.\kern-6pt
 \left[\vphantom{\copy0}\right.
 \kern-4pt
 \copy0
 \kern-4pt
 \left.\vphantom{\copy0}\right]\kern-6pt
 \left.\vphantom{\copy0}\right]_{#3}}

\newcommand{\jj}{\mathsf{j}}
\newcommand{\mm}{\mathsf{m}}

\newcommand{\rr}{\mathsf{r}}
\newcommand{\rrr}{\mathsf{r}}
\newcommand{\xx}{\mathsf{x}}
\newcommand{\ttt}{\mathsf{t}}
\newcommand{\zz}{\mathsf{z}}

\DeclareMathOperator{\ee}{e}

 \makeatletter
 \newcommand{\starit}[3]
 {\mkern #1mu {\vphantom{#3}}^*\mkern-#2 mu #3}
\renewcommand{\*}{\futurelet\@startemp\@starit}
\newcommand{\@starit}{%
 \def\starbefore{1} \def\starafter{6}
 \let\testit=\Iq \ifx \@startemp \testit
    \def\starbefore{0}\def\starafter{3} \fi
 \let\testit=\SL \ifx \@startemp \testit
    \def\starbefore{0}\def\starafter{1} \fi
 \let\testit=\Elem \ifx \@startemp \testit
    \def\starbefore{0}\def\starafter{4} \fi
 \let\testit=\Enorm \ifx \@startemp \testit
    \def\starbefore{0}\def\starafter{4} \fi
 \let\testit=\genset \ifx \@startemp \testit
    \def\starbefore{0}\def\starafter{5} \fi
 \let\testit=\natural \ifx \@startemp \testit
    \def\starbefore{0}\def\starafter{2} \fi
 \let\testit=B \ifx \@startemp \testit
    \def\starbefore{0}\def\starafter{3} \fi
 \let\testit=G \ifx \@startemp \testit
    \def\starbefore{0}\def\starafter{3} \fi
 \let\testit=X \ifx \@startemp \testit
    \def\starbefore{0}\def\starafter{4} \fi
 \let\testit=a \ifx \@startemp \testit
    \def\starbefore{0}\def\starafter{4} \fi
 \let\testit=\LU \ifx \@startemp \testit
    \def\starbefore{0}\def\starafter{5} \fi
 \let\testit=\BS \ifx \@startemp \testit
    \def\starbefore{0}\def\starafter{5} \fi
 \starit{\starbefore}{\starafter}}
  \makeatother

\newcommand{\pref}[1]{{\rm(}\ref{#1}{\rm)}}
\renewcommand{\see}[1]{{\rm(}see~\ref{#1}{\rm)}}
\newcommand{\seeand}[2]{{\rm(}see \ref{#1} and~\ref{#2}{\rm)}}
\newcommand{\seeSect}[1]{{\rm(}see \S\ref{#1}{\rm)}}
\newcommand{\seeDefn}[1]{{\rm(}see Definition~\ref{#1}{\rm)}}
\newcommand{\cf}[1]{{\rm(}cf.~\ref{#1}{\rm)}}
\newcommand{\cfand}[2]{{\rm(}cf.\ \ref{#1} and~\ref{#2}{\rm)}}
\newcommand{\cfSect}[1]{{\rm(}cf.\ \S\ref{#1}{\rm)}}
\newcommand{\fullref}[2]{\ref{#1}\pref{#1-#2}}
\newcommand{\fullsee}[2]{{\rm(}see~\fullref{#1}{#2}{\rm)}}

\newenvironment{thmref}{\thmrefer}{}
\newcommand{\thmrefer}[1]{\renewcommand\theequation
 {\protect\ref{#1}$'$}\addtocounter{equation}{-1}}

\numberwithin{equation}{section}

\makeatletter
\def\swappedhead@plain#1#2#3{%
 \textnormal{(\thmnumber{#2})}\thmname{ #1}\thmnote{ {\theoremnotefont#3}}}
 \newcommand{\theoremnotefont}{\normalfont}
 \makeatother
\swapnumbers

\newtheorem{thm}[equation]{Theorem}
\newtheorem{cor}[equation]{Corollary}
\newtheorem{prop}[equation]{Proposition}
\newtheorem{lem}[equation]{Lemma}
\newtheorem{Leibniz}[equation]{Leibniz' Principle}

\theoremstyle{definition}
 \newtheorem{defn}[equation]{Definition}
 \newtheorem{notation}[equation]{Notation}
 \newtheorem{assump}[equation]{Assumption}
 \newtheorem{rem}[equation]{Remark}
 \newtheorem{eg}[equation]{Example}

\newenvironment{ack}[1][\unskip]{
 \medskip \noindent \bf Acknowledgments. \rm}{\unskip\upshape\par}

\newenvironment{claim*}[1][\unskip]{
 \medskip \noindent \bf Claim. \it}{\unskip\upshape\par}

\newenvironment{case}[1][\unskip]{\refstepcounter{case}%
 \em
 \medskip \noindent Case \thecase\ #1.\ }{\unskip\upshape}
 \newcommand{\thecase}{\arabic{case}}
 \newcounter{case}

 \newcounter{step}
 \newenvironment{step}[1][\unskip]{\refstepcounter{step}
 \em
 \medskip \noindent Step \thestep\ #1.\ }{\unskip\upshape}
 \renewcommand{\thestep}{\arabic{step}}

\renewcommand{\thesubsection}{\thesection\Alph{subsection}}

\makeatletter
\def\subsection{\@startsection{subsection}{2}%
  \z@{-9pt plus -2.5pt minus -12pt}{-3pt}%
  {\normalfont\bfseries\S}}
 \makeatother

\makeatletter
 \renewcommand{\tocsubsection}[3]{%
  \indentlabel{\@ifnotempty{#2}{
 \hskip 0.125in \ignorespaces#1 \S#2.\enspace}}#3}
 \makeatother

\makeatletter
\newcommand{\@dotsep}{3}
\def\@tocline#1#2#3#4#5#6#7{\relax
  \ifnum #1>\c@tocdepth 
  \else
    \par \addpenalty\@secpenalty\addvspace{#2}%
    \begingroup \hyphenpenalty\@M
    \@ifempty{#4}{%
      \@tempdima\csname r@tocindent\number#1\endcsname\relax
    }{%
      \@tempdima#4\relax
    }%
    \parindent\z@ \leftskip#3\relax \advance\leftskip\@tempdima\relax
    \rightskip\@pnumwidth plus1em \parfillskip-\@pnumwidth
  #5\leavevmode\hskip-\@tempdima #6\relax
 \ifnum #1=1 
    \leaders\hbox{$\m@th
      \mkern \@dotsep mu\raise2pt\hbox{.}\mkern \@dotsep mu$}
 \else \ \fi\hfill
 $\mkern -\@dotsep mu$%
\hbox 
 {\@tocpagenum{#7}}\par
    \nobreak
    \endgroup
  \fi}

\begin{document}

\begin{abstract}
  We present unpublished work of D.~Carter, G.~Keller, and E.~Paige on bounded
generation in special linear groups. Let $n$~be a positive integer, and let $A
= \ints$ be the ring of integers of an algebraic number field~$K$ (or,
more generally, let $A$ be a localization $\localize{\ints}{S}$). If $n =
2$, assume that $A$~has infinitely many units.

We show there is a finite-index subgroup~$H$ of $\SL(n,A)$, such that every
matrix in~$H$ is a product of a bounded number of elementary matrices. We also
show that if $T \in \SL(n,A)$, and $T$ is not a scalar matrix, then there is a
finite-index, normal subgroup~$N$ of $\SL(n,A)$, such that every element
of~$N$ is a product of a bounded number of conjugates of~$T$.

For $n \ge 3$, these results remain valid when $\SL(n,A)$ is replaced by any of
its subgroups of finite index.
 \end{abstract}
 
\maketitle

\tableofcontents

 \section{Introduction} \label{0.}

This paper presents unpublished work of David Carter, Gordon Keller, and Eugene
Paige \cite{CKP} --- they should be given full credit for the results and the
methods of proof that appear here (but the current author is responsible for
errors and other defects in this manuscript). 
 Much of this work is at least 20 years old (note that it is mentioned in
\cite[p.~152 and bibliography]{DennisVaserstein}), but it has never been
superseded.

If a set~$\genset$ generates a group~$G$, then every element of~$G$ can be
written as a word in $\genset \cup \genset^{-1}$. We are interested in cases
where the length of the word can be bounded, independent of the particular
element of~$G$.

\begin{defn}
 A subset~$\genset$ of a group~$G$ \emph{boundedly generates}~$G$ if there is a
positive integer~$r$, such that every element of~$G$ can be written as a word
of length~$\le r$ in $\genset \cup \genset^{-1}$. 
  That is, for each $g \in G$, there is a sequence $x_1,x_2,\ldots,x_\ell$ of
elements of $\genset \cup \genset^{-1}$, with $\ell \le r$, such that 
 $ g = x_1 x_2 \cdots x_\ell $.
 \end{defn}

A well-known paper of D.~Carter and G.~Keller \cite{CarterKeller-BddElemGen}
proves that if $B$~is the ring of integers of a number field~$K$, and $n \ge
3$, then the set of elementary matrices $E_{i,j}(b)$ boundedly generates
$\SL(n,B)$. One of the two main results of \cite{CKP} is the following theorem
that generalizes this to the case $n = 2$, under an additional (necessary)
condition on~$B$. (For the proof, see
Corollary~\fullref{BddGen>3}{LU} and Theorem~\ref{BddGen2}.)

\begin{thm}[{(Carter-Keller-Paige \cite[(2.4) and (3.19)]{CKP})}]
\label{EijBddGen}
 Suppose 
 \begin{itemize}
 \item $B$ is the ring of integers of an algebraic number field~$K$ {\rm(}or,
more generally, $B$~is any order in the integers of~$K${\rm)},
 \item $n$~is a positive integer,
 \item $\Elem(n,B)$ is the subgroup of $\SL(n,B)$ generated by the elementary
matrices,
 and
 \item either $n \ge 3$, or $B$ has infinitely many units.
 \end{itemize}
 Then the elementary matrices boundedly generate $\Elem(n,B)$.

More precisely, there is a positive integer $r = r(n,k)$, depending only
on~$n$ and the degree~$k$ of~$K$ over~$\rational$, such that
 \begin{enumerate}
 \item every matrix in $\Elem(n,B)$ is a product of $\le r$ elementary matrices,
 and
 \item $\# \bigl( \SL(n,B)/\Elem(n,B) \bigr) \le r$.
 \end{enumerate}
 \end{thm}

\begin{rem} \ \label{EijBddGenRem}
 If $B$ is (an order in) the ring of integers of a number field~$K$, and $B$
has only finitely many units, then $K$ must be either $\rational$ or an
imaginary quadratic extension of~$\rational$. In this case, the elementary
matrices do not boundedly generate $\SL(2,B)$ \cite[Cor.\ of Prop.~8,
p.~126]{Tavgen-BddGenChev}.
 (This follows from the fact \cite{GrunewaldSchwermer} that some finite-index
subgroup of $\SL(2,B)$ has a nonabelian free quotient.) Thus, our assumption
that $n \ge 3$ in this case is a necessary one.
 \end{rem}

The following result is of interest even when
$\genset$ consists of only a single matrix~$X$.

\begin{thmref}{BddGenNormal}
 \begin{thm}[{(Carter-Keller-Paige \cite[(2.7) and (3.21)]{CKP})}]
\label{ConjBddGen}
 Let
 \begin{itemize}
 \item $B$ and $n$ be as in Theorem~\ref{EijBddGen},
 \item $\genset$ be any subset of $\SL(n,B)$ that does not consist
entirely of scalar matrices,
 and
 \item $\genset^{\normal} = \bigset{ T^{-1} X T }{ 
 \begin{matrix} X \in \genset, \\ T \in \SL(n,B) \end{matrix} }$.
 \end{itemize}
 Then $\genset^{\normal}$ boundedly generates a finite-index normal
subgroup of\/ $\SL(n,B)$.
 \end{thm}
 \end{thmref}

\begin{rem} \ 
 \begin{enumerate}
 \item In the situation of Theorem~\ref{ConjBddGen}, let
$\gen{\genset^{\normal}}$ be the subgroup generated by~$\genset^{\normal}$. It
is obvious that $\genset^{\normal}$ is a normal subgroup of $\SL(n,B)$, and it
is well known that this implies that $\genset^{\normal}$ has finite index in
$\SL(n,B)$ (cf.\ \ref{Sandwich}, \ref{u4inN}, and \ref{SLn/EnormFinite}).
 \item The conclusion of Theorem~\ref{ConjBddGen} states that there is a positive
integer~$r$, such
that every element of $\genset^{\normal}$ is a product of $\le r$ elements of
$\genset^\normal$ (and their inverses).
 Unlike in \pref{EijBddGen}, we do \emph{not} prove that the bound~$r$ can be
chosen to depend on only $n$ and~$k$. See Remark~\ref{normalbound} for a
discussion of this issue.
 \item We prove Thms.~\ref{EijBddGen} and \ref{ConjBddGen} in a more general
form that allows $B$ to be replaced with any localization $\BS$. 
 It is stated in \cite{CKP} (without proof) that the same conclusions hold if
$B$ is replaced by an arbitrary subring~$A$ of any number field (with the
restriction that $A$~is required to have infinitely many units if $n = 2$). It
would be of interest to establish this generalization.
 \item If $\Gamma$ is any subgroup of finite index in $\SL(n,B)$, then
Theorem~\fullref{BddGenNormal}{Gamma} is a generalization of Theorem~\ref{EijBddGen}
that applies with $\Gamma$ in the place of $\SL(n,B)$.
 For $n \ge 3$, Theorem~\ref{BddGenGammaN} is a generalization of
Theorem~\ref{ConjBddGen} that applies with $\Gamma$ in the place of $\SL(n,B)$.
 \end{enumerate}
 \end{rem}

Let us briefly outline the proof of Theorem~\ref{EijBddGen}. (A similar approach
applies to Theorem~\ref{ConjBddGen}.) For $n$ and~$B$ as in the statement of the
theorem, it is known that the subgroup $\Elem(n,B)$ generated by the elementary
matrices has finite index in $\SL(n,B)$ \cite{BMS, Serre-CSPSL2,
Vaserstein-SL2}. Theorem~\ref{EijBddGen} is obtained by axiomatizing this
proof: 
 \begin{enumerate}
 \item Certain ring-theoretic axioms are defined (for $n \ge 3$, the axioms are
called $\SRquot$, $\GEN(\ttt,\rr)$, and $\EXP(\ttt,\ell)$, where the
parameters $\ttt$, $\rr$, and~$\ell$ are positive integers).
 \item It is shown that the ring~$B$ satisfies these axioms (for appropriate
choices of the parameters).
 \item It is shown that if $A$ is any integral domain satisfying these axioms,
then $\Elem(n,A)$ is a finite-index subgroup of $\SL(n,A)$.
 \end{enumerate}
 The desired conclusion is then immediate from the following simple
consequence of the Compactness Theorem of first-order logic
\seeSect{CpctnessSect}:

\begin{prop} \label{Cpct:E(n,A;q)}
 Let
 \begin{itemize}
 \item $n$ be a positive integer,
 and
 \item $\theory$ be a set of first-order axioms in the language of ring
theory.
 \end{itemize}
 Suppose that, for every commutative ring~$A$ satisfying the axioms
in~$\theory$, the subgroup $\Elem(n,A)$ generated by the elementary matrices has
finite index in $\SL(n,A)$.
 Then, for all such~$A$, the elementary matrices boundedly generate $\Elem(n,A)$.

More precisely, there is a positive
integer $r = r(n,\theory)$, such that, for all~$A$ as above, every matrix in
$\Elem(n,A)$ is a product of $\le r$ elementary matrices.
 \end{prop}

\begin{eg}
 It is a basic fact of linear algebra that if $F$ is any field, then every
element of $\SL(n,F)$ is a product of elementary matrices.
 This yields the conclusion that $\Elem(n,F) = \SL(n,F)$. Since fields are
precisely the commutative rings satisfying the additional axiom $(\forall x)
(\exists y) (x \neq 0 \rightarrow xy = 1)$, then Proposition~\ref{Cpct:E(n,A;q)}
implies that each element of $\SL(n,F)$ is the product of a bounded number of
elementary matrices. (Furthermore, a bound on the number of elementary
matrices can be found that depends only on~$n$, and is universal for all
fields.) In the case of fields, this can easily be proved directly, by
counting the elementary matrices used in a proof that $\Elem(n,F) = \SL(n,F)$, but
the point is that this additional work is not necessary --- bounded generation
is an automatic consequence of the fact that $\Elem(n,A)$ is a finite-index
subgroup. 
 \end{eg}

Because we obtain bounded generation from the Compactness Theorem (as in
\pref{Cpct:E(n,A;q)}), the conclusions in this paper do not provide any
explicit bounds on the number of matrices needed. It should be possible to
obtain an explicit formula by carefully tracing through the arguments in this
paper and in the results that are quoted from other sources, but this would
be nontrivial (and would make the proofs messier).
The applications we have in mind do not require this.

\begin{rem}
 Assuming a certain strengthening of the Riemann Hypothesis, Cooke and
Weinberger \cite{CookeWeinberger} proved a stronger version of
Theorem~\ref{EijBddGen} that includes an explicit estimate on the integer~$r$
(depending only on~$n$, not on~$k$), under the assumption that $B$ is the full
ring of integers, not an order. 
 For $n \ge 3$, the above-mentioned work of D.~Carter and G.~Keller
\cite{CarterKeller-BddElemGen, CarterKeller-ElemExp} removed the reliance on
unproved hypotheses, but obtained a weaker bound that depends on the
discriminant of the number field.
 For $n=2$, B.~Liehl \cite{Liehl-Beschrankte} proved bounded generation (without explicit
bounds), but required some assumptions on the number field~$K$.
 More recently, for a localization $B_{\mathcal{S}}$ with $\mathcal{S}$ a
sufficiently large set of primes, D.~Loukanidis and V.~K.~Murty
\cite{LoukanidisMurty, Murty-BddGen} obtained explicit bounds for
$\SL(n,B_{\mathcal{S}})$ that depend only on $n$ and~$k$, not the
discriminant. 

There is also interesting literature on bounded generation of other
(arithmetic) groups, e.g., \cite{AdianMennicke, Bardakov, DennisVaserstein,
ErovenkoRapinchukCR, ErovenkoRapinchukPreprint, LoukanidisMurty, Murty-BddGen,
Rapinchuk-CSPFinWidth, Shalom-BddGen, SivatskiStepanov, Tavgen-BddGenChev, 
Tavgen-FinWidthArith, vanderKallen-SL3(C[X]), Zakiryanov}.
 \end{rem}


\begin{ack}
 This paper was written during a visit to the University of Auckland. I would
like to thank the Department of Mathematics of that institution for its
hospitality. I would also like to thank Jason Manning, Lucy Lifschitz, and
Alex Lubotzky for bringing the preprint \cite{CKP} to my attention,
and an anonymous referee for reading the manuscript carefully and
providing numerous corrections and helpful comments.
 The work was partially supported by a grant from the National Sciences and
Engineering Research Council of Canada.
 \end{ack}

\section{Preliminaries}

\begin{assump}
 All rings are assumed to have~$1$, and any subring is assumed to contain the
multiplicative identity element of the base ring. (This is taken to be part of
the definition of a ring or subring.)
 \end{assump}

\subsection{Notation} \label{0.3.}

\begin{defn}
 Let $B$ be an integral domain.
 \begin{enumerate}
 \item A subset~$S$ of~$B$ is  \emph{multiplicative} if
 $S$ is closed under multiplication,
 and
 $0 \notin S$.
 \item If $S$ is a multiplicative subset of~$B$, then 
 $$ \BS = \bigset{ \frac{b}{s} }{ b \in B, s \in S} .$$
 This is a subring of the quotient field of~$B$.
 \end{enumerate}
 \end{defn}

As usual, we use $\gen{\genset}$ to denote the subgroup generated by
a subset~$\genset$ of a group~$G$.
 In order to conveniently discuss bounded generation, we augment this
notation with a subscript, as follows. 

\begin{defn}
 For any subset $\genset$ of a group~$G$, and any nonnegative integer~$r$, we
define $\gennum{\genset}{r}$, inductively, by:
 \begin{itemize}
 \item $\gennum{\genset}{0} = \{1\}$ (the identity element of~$G$), 
 and
 \item $\gennum{\genset}{r+1} = \gennum{\genset}{r} \cdot \left( \genset \cup
\genset^{-1} \cup \{1\} \right)$.
 \end{itemize}
 That is, $\gennum{\genset}{r}$ is the set of elements of~$G$ that can be
written as a word of length $\le r$ in $\genset \cup \genset^{-1}$.
 Thus, $\genset$ boundedly generates~$G$ if and only if we have
$\gennum{\genset}{r} = G$, for some positive integer~$r$.
 \end{defn}

\begin{notation} \label{PrelimNotn}
 Let $A$ be a commutative ring,
 $\Iq$ be an ideal of~$A$,
 and
 $n$ be a positive integer.
 \begin{enumerate}
 \itemsep=\smallskipamount 

 \item $\Id{n}$ denotes the $n \times n$ identity matrix.

 \item $\SL(n,A;\Iq) = \{\, T \in \SL(n,A) \mid T \equiv \Id{n} \mod \Iq \,\}$.

 \item For $a \in A$, and $1 \le i,j \le n$ with $i \neq j$, we use
$E_{i,j}(a)$ to denote the $n \times n$ \emph{elementary matrix}, such that the
only nonzero entry of $E_{i,j}(a) - \Id{n}$ is the $(i,j)$~entry, which is~$a$.
(We may use $E_{i,j}$ to denote $E_{i,j}(1)$.)

 \item $\LU(n,\Iq) = \bigset{ E_{i,j}(a) }{
 \begin{matrix} a \in \Iq, \\ 1 \le i,j \le n, \\ i \neq j \end{matrix}
 }$.
 In other words, $\LU(n,A)$ is the set of all $n \times n$ elementary
matrices, and $\LU(n,\Iq) = \LU(n,A) \cap \SL(n,A;\Iq)$.

 \item $\Elem(n,\Iq) = \langle \LU(n,\Iq) \rangle$.
 Thus, $\Elem(n,A)$ is the subgroup of $\SL(n,A)$ generated by the elementary
matrices.

 \item $\LUnorm(n,A;\Iq)$ is the set of $\Elem(n,A)$-conjugates of elements of
$\LU(n,\Iq)$.

 \item $\Enorm(n,A;\Iq) = \langle \LUnorm(n,A;\Iq) \rangle$.
 Thus, $\Enorm(n,A;\Iq)$ is the smallest \emph{normal} subgroup of $\Elem(n,A)$
that contains $\LU(n,\Iq)$.

\item  \label{PrelimNotn-W(q)}
 $\displaystyle W(\Iq) = \bigset{ (a,b) \in A \times A }{
 \begin{matrix}
 (a,b) \equiv (1,0) \mod \Iq
 \\ \text{and} \\
 a A + bA = A \end{matrix}
 }$.
 Note that $(a,b) \in W(\Iq)$ if and only if there exist $c,d \in A$, such
that 
 $\mat{a}{b}{c}{d} \in \SL(n,A;\Iq)$ \cite[Prop.~1.2(a),
p.~283]{Bass-AlgKTheoryBook}.

 \item $\U{\Iq}$ is the group of units of~$A /
\Iq$.

 \end{enumerate}
 \end{notation}

 Note that $\Elem(n,A)$ is boundedly generated by elementary
matrices if and only if $\Elem(n,A) = \gennum{\LU(n,A)}{r}$, for some positive
integer~$r$.

\begin{rem}
 The subgroup $\Enorm(n,A;\Iq)$ is usually denoted $\Elem(n,A;\Iq)$
in the literature, but we include the superscript ``$\normal$'' to emphasize
that this subgroup is normalized by $\Elem(n,A)$, and thereby reduce
the likelihood of confusion with $\Elem(n,\Iq)$.
 \end{rem} 

\begin{notation}
 Suppose $K$ is an algebraic number field. We use $\Norm =
\Norm_{K/\rational}$ to denote the norm map from~$K$ to~$\rational$.
 \end{notation}

\subsection{The Compactness Theorem of first-order logic} \label{CpctnessSect}

The well-known G\"odel Completeness Theorem states that if a theory in
first-order logic is consistent (that is, if it does not lead to a
contradiction of the form $\varphi \wedge \neg \varphi$), then the theory
has a model. Because any proof must have finite length, it can quote only
finitely many axioms of the theory. This reasoning leads to the following
fundamental theorem, which can be found in introductory texts on
first-order logic.

\begin{thm}[(Compactness Theorem)] \label{Cpctness}
 Suppose $\theory$ is any set of first-order sentences {\rm(}with no free
variables{\rm)} in some first-order language~$\lang$.
 If $\theory$ does not have a model, then some finite subset~$\theory_0$
of~$\theory$ does not have a model.
 \end{thm}

\begin{cor} \label{CpctnessCor}
 Fix a positive integer~$n$, and let $\lang$ be a first-order language that
contains
 \begin{itemize}
 \item the language of rings {\rm(}$+, \times, 0, 1${\rm)},
 \item $n^2$ variables $x_{ij}$ for $1 \le i,j \le n$,
 \item two $n^2$-ary relation symbols $X(x_{ij})$ and $H(x_{ij})$,
 and
 \item any number {\rm(}perhaps infinite{\rm)} of other variables, constant
symbols, and relation symbols.
 \end{itemize}
 Suppose $\theory$ is a set of sentences in the language~$\lang$, such that,
for every model 
 $$\bigl(A,(+,\times,0,1,X, H,\ldots )\bigr)$$
 of the theory~$\theory$,
 \begin{itemize}
 \item the universe $A$ is a commutative ring {\rm(}under the binary operations
$+$ and~$\times${\rm)},
 and
 \item letting 
 $$ \text{$\displaystyle
 X_A = \bigset{(a_{ij})_{i,j=1}^n }
 { \begin{matrix} a_{ij} \in A, \\ X(a_{ij}) \end{matrix} }$
 \quad and \quad
 $\displaystyle
 H_A = \bigset{(a_{ij})_{i,j=1}^n }
 { \begin{matrix} a_{ij} \in A, \\ H(a_{ij}) \end{matrix} }$,}$$
 we have
 \begin{itemize} 
 \item $H_A$ is a subgroup of\/ $\SL(n,A)$,
 and
 \item $X_A$ generates a subgroup of finite index in~$H_A$.
 \end{itemize}
 \end{itemize}
 Then, for every model $\bigl(A, \ldots\bigr)$ of~$\theory$, the set $X_A$
\textbf{boundedly generates} a subgroup of finite index in~$H_A$.

More precisely, there is a positive integer $r = r(n,\lang,\theory)$, such
that, for every model $(A,\ldots)$ of~$\theory$, 
 $\gennum{X_A}{r}$ is a subgroup of~$H_A$, and the index of this subgroup
is~$\le r$.
 \end{cor}
 
 \begin{proof}
 This is a standard argument, so we provide only an informal sketch.
 \begin{itemize}
 \item Let $\lang^+$ be obtained from~$\lang$ by adding constant symbols to
represent infinitely many matrices $C_1,C_2,C_3,\ldots$. (Each matrix requires 
$n^2$~constant symbols $c_{i,j}$.)
 \item Let $\theory^+$ be obtained from~$\theory$ by adding first-order
sentences specifying, for all $i,j,r \in \natural^+$, with $i \neq j$, that
	\begin{itemize}
	\item $C_i \in H_A$, and
	\item $C_i^{-1} C_j \notin \gennum{X_A}{r-1}$.
	\end{itemize}
 \end{itemize}
Since $X_A$ generates a subgroup of finite index in~$H_A$, we know that
$\theory^+$ is not consistent. From the Compactness Theorem, we conclude,
for some~$r$, that it is impossible to find $C_1,C_2,\ldots,C_r \in H_A$, such
that $C_i^{-1} C_j \notin \gennum{X_A}{r-1}$ for $i \neq j$. This implies the
index of $\langle X_A \rangle$ is less than~$r$. Also, we must have
$\gennum{X_A}{r^2} = \langle X_A \rangle$ (otherwise, we could choose
$C_i \in \gennum{X_A}{i r} \smallsetminus \gennum{X_A}{i r - 1}$).
 \end{proof}

\begin{proof}[Proof of Proposition~\ref{Cpct:E(n,A;q)}]
 This is a standard compactness argument, so we provide only a sketch. Let
$\theory'$ consist of:
 \begin{itemize}
 \item the axioms in~$\theory$,
 \item the axioms of commutative rings,
 \item a collection of sentences that guarantees $X_A = \LU(n,A)$,
 and
 \item  a collection of sentences that guarantees $H_A = \SL(n,A)$.
 \end{itemize}
 Then the desired conclusion is immediate from Corollary~\ref{CpctnessCor}.
 \end{proof}

\subsection{Stable range condition $\SR_{\mm}$}

We recall the stable range condition $\SR_{\mm}$ of Bass. (We use the
indexing of \cite{HahnOMeara}, not that of \cite{Bass-AlgKTheoryBook}.) For
convenience, we also introduce a condition $\SRquot$ that is intermediate
between $\SR_1$ and $\SR_2$. In our applications, the parameter $\mm$ will
always be either $1$ or $1\frac{1}{2}$ or~$2$.

\begin{defn}[{(\cite[Defn.~3.1, p.~231]{Bass-AlgKTheoryBook},
\cite[p.~142]{HahnOMeara}, cf.\ \cite[\S4]{Bass-KTheoryStable})}]
\label{SRDefn}
 Fix a positive integer~$\mm$. 
 We say that a commutative ring~$A$ satisfies the \emph{stable range} condition
$\SR_\mm$ if,
 for all $a_0,a_1,\ldots,a_r \in A$, such that
 \begin{itemize}
 \item $r \ge \mm$
 and 
 \item $a_0 A + a_1 A +\cdots + a_r A = A$,
 \end{itemize}
 there exist $a_1',a_2',\ldots,a_{r}' \in A$, such that 
 \begin{itemize}
 \item $a_i' \equiv a_i \mod a_0 A$, for $1 \le i \le r$,
 and
 \item $a_1' A + \cdots + a_{r}' A = A$.
 \end{itemize}
 \end{defn}

The condition $\SR_\mm$ can obviously be represented by a list of infinitely
many first-order statements, one for each integer $r \ge \mm$. It is
interesting (though not necessary) to note that the single case $r = \mm$
implies all the others \cite[(4.1.7), p.~143]{HahnOMeara}, so a single
statement suffices.

\begin{defn}
 We say a commutative ring~$A$ satisfies $\SRquot$ if $A/\Iq$ satisfies
$\SR_1$, for every nonzero ideal~$\Iq$ of~$A$.
 \end{defn}

It is easy to see that $\SR_1 \Rightarrow \SRquot \Rightarrow \SR_2$.

\begin{rem}
 If $A$ satisfies $\SR_\mm$ (for some~$\mm$), and $\Iq$ is any ideal of~$A$,
then $A/\Iq$ also satisfies $\SR_\mm$ \cite[Lem.~4.1]{Bass-KTheoryStable}.
Hence, $A$ satisfies $\SRquot$ if and only if $A/qA$ satisfies $\SR_1$, for
every nonzero $q \in A$. This implies that $\SRquot$ can be expressed in terms
of first-order sentences.
 \end{rem}

\begin{notation}
 As is usual in this paper,
 \begin{itemize}
 \item $K$ is an algebraic number field,
 \item $\ints$ is the ring of integers of~$K$, 
 \item $B$~is an order in~$\ints$,
 and
 \item $S$~is a multiplicative subset of~$B$.
 \end{itemize}
 \end{notation}

The following result is well known.

\begin{lem}
 $\BS$ satisfies $\SRquot$.
 \end{lem}

\begin{proof}
 Let $\Iq$ be any nonzero ideal of~$\BS$. Since the quotient ring $\BS/\Iq$ is
finite, it is semilocal. So it is easy to see that it satisfies $\SR_1$
\cite[Prop.~2.8]{Bass-AlgKTheoryBook}.
 \end{proof}

The following fundamental result of Bass is the reason for our interest in
$\SR_\mm$.

\begin{thm}[{(Bass \cite[\S4]{Bass-KTheoryStable})}] \label{BassThm}
 Let
 \begin{itemize}
 \item $A$ be a commutative ring,
 \item $\mm$ be a positive integer, such that $A$ satisfies the stable range
condition $\SR_\mm$,
 \item $n > \mm$,
 and
 \item $\Iq$ be an ideal of~$A$.
 \end{itemize}
 Then:
 \begin{enumerate}
 \item \label{BassThm-GL=GLxE}
 $\SL(n,A;\Iq) = \SL(\mm,A;\Iq) \Enorm(n,A;\Iq)$.
 \item $\Enorm(n,A;\Iq)$ is a normal subgroup of\/ $\SL(n,A)$.
 \item If $n \ge 3$, then $\bigl[ \Elem(n,A), \SL(n,A;\Iq) \bigr] =
\Enorm(n,A;\Iq)$.
 \end{enumerate}
 \end{thm}

Applying the case $\mm = 1$ of \fullref{BassThm}{GL=GLxE} to the quotient ring
$A/\Iq'$ yields the following conclusion:

\begin{cor} \label{VasLem1Normal}
 Let
 \begin{itemize}
 \item $A$ be a commutative ring,
 \item $n$ be a positive integer,
 and
 \item $\Iq$ and~$\Iq'$ be nonzero ideals of~$A$, such that $\Iq' \subseteq
\Iq$.
 \end{itemize}
 If $A/\Iq'$ satisfies $\SR_1$, then $\SL(n,A;\Iq) = \SL(n,A;\Iq')\Enorm(n,A;
\Iq)$.
 \end{cor}

\subsection{Mennicke symbols}

We recall the definition and basic properties of Mennicke symbols, including
their important role in the study of the quotient group $\SL(n,A;\Iq)/\Enorm(n,A;\Iq)$.

 \begin{defn}[{\cite[Defn.~2.5]{BMS}}] \label{1.1}
 Suppose $A$ is a commutative ring and $\Iq$ is an ideal in~$A$. 
 Recall that $W(\Iq)$ was defined in \fullref{PrelimNotn}{W(q)}.
 \begin{enumerate}
 \item A \emph{Mennicke symbol} is a function $(a,b) \mapsto \Men{b}{a}{}$
from $W(\Iq)$ to a group~$C$, such that
 \newcounter{firstMS} \newcounter{secondMS}
 \renewcommand{\theequation}{MS\arabic{firstMS}\alph{secondMS}}
 \begin{align}
 \setcounter{firstMS}{1} \setcounter{secondMS}{1}
 \label{MS1a} 
 \Men{b+ta}{a}{} &= \Men{b}{a}{}
 && \text{whenever $(a,b) \in W(\Iq)$ and $t \in \Iq$;} \\
 \setcounter{firstMS}{1} \setcounter{secondMS}{2}
 \label{MS1b}
 \Men{b}{a+tb}{} &= \Men{b}{a}{}
 &&\text{whenever $(a,b) \in W(\Iq)$ and $t \in A$;
 and} \\
 \setcounter{firstMS}{2} \setcounter{secondMS}{1}
 \label{MS2a}
 \Men{b_1}{a}{} \Men{b_2}{a}{} &= \Men{b_1 b_2}{a}{}
 &&\text{whenever $(a,b_1), (a,b_2) \in W(\Iq)$.}
 \end{align}
 \addtocounter{equation}{-3} 
 \item It is easy to see that, for some group~$C(\Iq)$ (called the
\emph{universal Mennicke group}), there is a \emph{universal} Mennicke symbol
 $$ \UMen{\ }{\ }{\Iq} \colon W(\Iq) \to C(\Iq) ,$$
 such that any Mennicke symbol $\Men{\ }{\ }{} \colon W(\Iq) \to C$, for
any group~$C$, can be obtained by composing $\UMen{\ }{\ }{\Iq}$ with a unique
homomorphism from~$C(\Iq)$ to~$C$. The universal Mennicke symbol and the
universal Mennicke group are unique up to isomorphism.
 \end{enumerate}
 \end{defn}

 The following classical theorem introduces Mennicke symbols into the
study of $\Enorm(n,A;\Iq)$.

\begin{notation}
 For convenience, when $T \in \SL(2,A;\Iq)$, we use $\overline{T}$ to denote
the image of~$T$ under the usual embedding of $\SL(2,A;\Iq)$ in the top left
corner of $\SL(n,A;\Iq)$.
 \end{notation}

\begin{thm}[{\cite[Thm.~5.4 and Lem.~5.5]{BMS},
\cite[Prop.~1.2(b), p.~283 and Thm.~2.1(b), p.~293]{Bass-AlgKTheoryBook}}]
\label{MennickeThm}
 Let
 \begin{itemize}
 \item $A$ be a commutative ring,
 \item $\Iq$ be an ideal of~$A$,
 \item $N$ be a normal subgroup of $\SL(n,A;\Iq)$, for some $n \ge 2$,
 and
 \item $C = \SL(n,A;\Iq)/N$, \end{itemize}
 such that $N$ contains both $\Enorm(n,A;\Iq)$ and $\bigl[ \Elem(n,A),
\SL(n,A;\Iq) \bigr]$.
 Then:
 \begin{enumerate}
 \item The map $\Men{b}{a}{\Iq} \colon W(\Iq) \to C$, defined by
 $$ (a,b) \mapsto \Men{b}{a}{\Iq} = \overline{\mat{a}{b}{*}{*}} N ,$$
 is well-defined.
 \item $\Men{\ }{\ }{\Iq}$ satisfies \pref{MS1a} and \pref{MS1b}.
 \item {\rm (Mennicke)} \label{MennickeThm-MS2a}
 If $n \ge 3$, then $\Men{\ }{\ }{\Iq}$ also satisfies \pref{MS2a}, so it is a
Mennicke symbol.
 \end{enumerate}
 \end{thm}

Under the assumption that $A$ is a Dedekind ring,  Bass, Milnor, and Serre
\cite[\S2]{BMS} proved several basic properties of Mennicke symbols; these
results appear in \cite{Bass-AlgKTheoryBook} with the slightly weaker
hypothesis that $A$ is a Noetherian ring of dimension $\le 1$. For our
applications, it is important to observe that the arguments of
\cite{Bass-AlgKTheoryBook} require only the assumption that $A/\Iq$ satisfies
the stable range condition $\SR_1$, for every nonzero ideal~$\Iq$ of~$A$.

\begin{lem}[{(cf.\ \cite[\S2]{BMS}, \cite[\S6.1]{Bass-AlgKTheoryBook})}]
\label{BasicMennicke}
 Suppose
 \begin{itemize}
 \item $A$ is an integral domain that satisfies $\SRquot$,
 \item $\Iq$ is an ideal in~$A$,
 and
 \item $\Men{\ }{\ }{} \colon W(\Iq) \to C$ is a Mennicke symbol.
 \end{itemize}
 Then:
 \begin{enumerate}

 \item \label{BasicMennicke-10}
 $\Men{0}{1}{} = 1$ {\rm(}the identity element of~$C${\rm)}.

 \item \label{BasicMennicke-b(1-a)}
 If $(a,b) \in W(\Iq)$, then $\Men{b}{a}{} = \Men{b(1-a)}{a}{}$.

 \item \label{BasicMennicke-unit}
 If\/ $(a,b) \in W(\Iq)$, and there is a unit $u \in A$, such that either $a
\equiv u \mod bA$ or $b \equiv u \mod aA$, then $\Men{b}{a}{}
= 1$.

 \item \label{BasicMennicke-smaller}
 If\/ $(a,b) \in W(\Iq)$, and $\Iq'$ is any nonzero ideal contained in~$\Iq$,
then there exists $(a',b') \in W(\Iq')$, such that
$\displaystyle\Men{b}{a}{} = \Men{b'}{a'}{}$.



 \item \label{BasicMennicke-abelian}
  The image of the Mennicke symbol\/ $\Men{\ }{\ }{}$ is an abelian subgroup
of~$C$.


 \item {\rm (Lam)} \label{BasicMennicke-Lam}
 If\/ $\Iq$ is principal, then 
 \begingroup 
 \renewcommand{\theequation}{MS2b}
 \begin{equation}  \label{MS2b}
 \text{$\Men{b}{a_1}{} \Men{b}{a_2}{} = \Men{b}{a_1 a_2}{}$ whenever $(a_1,b),
(a_2,b) \in W(\Iq)$.}
 \end{equation}
 \endgroup 
 \addtocounter{equation}{-1} 

 \end{enumerate}
 \end{lem}

 The following result provides a converse to
Lemma~\fullref{BasicMennicke}{Lam}. It will be used in the proof of
Lemma~\ref{MenSL2CanPrinc}. 
 

\begin{lem} \label{MS2b->MS2a}
 Suppose 
 \begin{itemize}
 \item $A$ is a commutative ring,
 \item $\Iq$ is an ideal in~$A$,
 \item $C$ is a group,
 and
 \item $\Men{\ }{\ }{} \colon W(\Iq) \to C$ satisfies \pref{MS1a} and
\pref{MS1b}.
 \end{itemize}
 Then:
 \begin{enumerate}
 \item \label{MS2b->MS2a-full}
 {\rm (Lam \cite[Prop.~1.7(a), p.~289]{Bass-AlgKTheoryBook})} If\/ $\Men{\ }{\
}{}$ satisfies \pref{MS2b}, then it also satisfies \pref{MS2a}, so it is a
Mennicke symbol.
 \item \label{MS2b->MS2a-triv}
 If\/ $\Men{\ }{\ }{}$ satisfies \pref{MS2b} whenever $\Men{b}{a_2}{} = 1$,
then it satisfies \pref{MS2a} whenever $\Men{b_2}{a}{} = 1$.
 \end{enumerate}
 \end{lem}

\begin{proof}
 \pref{MS2b->MS2a-full} Given $\Men{b_1}{a}{}, \Men{b_2}{a}{} \in W(\Iq)$, let
$q = 1 - a \in \Iq$.
 Note that, for any $b \in \Iq$, we have
 \begin{equation} \label{MS2b->MS2aPf-btn}
 \text{$ \Men{b q^n}{a}{} = \Men{b}{a}{}$ for every positive integer~$n$}
 \end{equation} (because the proof of \fullref{BasicMennicke}{b(1-a)} does not
appeal to \pref{MS2a}). Also, because
 $$ \Men{bq^2}{1 + bq}{}
 = \Men{bq^2 - q(1 + bq)}{1 + bq}{}
 = \Men{-q}{1 + bq}{}
 = \Men{-q}{1}{}
 = 1 ,$$
 we have
 \begin{align} \label{MS2b->MS2aPf-recip}
 \Men{bq^2}{a}{}
 &= \Men{bq^2}{a}{} \Men{bq^2}{1 + bq}{} 
 = \Men{bq^2}{a(1 + bq)}{} 
 = \Men{bq^2}{a + abq}{} 
 \notag\\&= \Men{bq^2}{a + abq - bq^2}{} 
 = \Men{bq^2}{a + bq(a-q)}{} 
 \\&= \Men{bq^2}{a + bq(1)}{} 
 = \Men{bq^2 - q(a+bq)}{a + bq}{} 
 = \Men{-aq}{a + bq}{} 
 \notag. \end{align}
 Applying, in order,
 \pref{MS2b->MS2aPf-btn} to both factors, 
 \pref{MS2b->MS2aPf-recip} to both factors,
 \pref{MS2b},
 \pref{MS1b},
 definition of~$q$,
 \pref{MS1b},
 \pref{MS2b->MS2aPf-recip},
 and
 \pref{MS2b->MS2aPf-btn}, 
 yields
 \begin{align*}
 \Men{b_1}{a}{} \Men{b_2}{a}{}
 &= \Men{b_1 q^2}{a}{} \Men{b_2 q^2}{a}{}
 = \Men{-aq}{a + b_1 q}{} \Men{-aq}{a + b_2 q}{}
 \\&= \Men{-aq}{(a + b_1 q)(a + b_2 q)}{}
 = \Men{-aq}{a^2 + b_1 b_2 q^2}{}
 \\&= \Men{-aq}{a(1-q) + b_1 b_2 q^2}{}
 = \Men{-aq}{a + b_1 b_2 q^2}{}
 \\&= \Men{b_1 b_2 q^3}{a}{}
 = \Men{b_1 b_2}{a}{}
 . \end{align*}

\pref{MS2b->MS2a-triv}
 The condition \pref{MS2b} was applied only twice in the proof
of~\pref{MS2b->MS2a-full}.
 \begin{itemize}
 \item In the first application, the second factor is $\Men{bq^2}{1 + bq}{} =
1$.
 \item In the other application, the second factor is $\Men{-a q}{a + b_2 q}{}
= \Men{b_2}{a}{}$, which is assumed to be~$1$.
 \end{itemize}
 Therefore, exactly the same calculations apply. 
 \end{proof}

The following useful result is stated with a slightly weaker
hypothesis in \cite{Bass-AlgKTheoryBook}:

\begin{prop}[{\cite[Thm.~VI.2.1a, p.~293]{Bass-AlgKTheoryBook}}]
\label{MS2aTriv->MS2a}
 If $\Men{\ }{\ }{}$ satisfies \pref{MS2a} whenever $\Men{b_2}{a}{} = 1$, then
it is a Mennicke symbol.
 \end{prop}

Combining this with Lemma~\fullref{MS2b->MS2a}{triv} yields the following
conclusion:

\begin{cor} \label{MS2bTriv->MS2a}
 If $\Men{\ }{\ }{}$ satisfies \pref{MS2b} whenever $\Men{b}{a_2}{} = 1$, then
it is a Mennicke symbol.
 \end{cor}

We conclude this discussion with two additional properties of Mennicke symbols.

\begin{lem} \label{AddlMen}
 Let $A$, $\Iq$, and $\Men{\ }{\ }{}$ be as in Lemma~\ref{BasicMennicke}.
 \begin{enumerate}
 \item  \label{AddlMen-inverse}
 If $\mat{a}{b}{c}{d} \in \SL(2,A;\Iq)$, then $\Men{b}{a}{}^{-1} =
\Men{c}{a}{}$.
 \item \label{AddlMen-1.4} 
 Suppose
 $\Iq = qA$ is principal, and
 $a$, $b$, $c$, $d$, $f$,  and~$g$ are elements of~$A$, such that
 $$ \text{$\mat{a}{b}{c}{d}$ and $f \Id{2} + g \mat{a}{b}{c}{d}$ are in
$\SL(2,A; qA)$.} $$
 Then
 $$ \Men{bg}{f+ga}{}^2 = \Men{b}{f+ga}{}^2
.$$
 \end{enumerate}
 \end{lem}

\begin{proof}

 \pref{AddlMen-inverse} We have
 \begin{align*}
 \Men{b}{a}{} \Men{c}{a}{}
 &= \Men{bc}{a}{}
 = \Men{bc(1-a)}{a}{}
 \\&= \Men{(bc - ad)(1-a)}{a}{}
 = \Men{-(1-a)}{a}{}
 = \Men{a-1}{1}{}
 = 1 
 . \end{align*}

 \pref{AddlMen-1.4}
 Note that, by assumption, $a$, $d$, and $f + ga$ are all congruent to~$1$
modulo $qA$. Also, working modulo~$g q A$, we have
 \begin{align*}
 (f+ga)^2 
 &\equiv (f+g)^2 
 &&\text{(since $a \equiv 1  \mod q A$)}
 \\&\equiv f^2 + (a + d) fg + g^2
 &&\text{(since $a + d \equiv 1 + 1 = 2 \mod q A$)}
 \\&= \det \left( f \Id{2} + g \mat{a}{b}{c}{d} \right)
 \\&= 1
 . \end{align*}
 Therefore
 \begin{align*}
 \Men{bg}{f+ga}{}^2
 &= \Men{bg}{(f+ga)^2}{}
 && \text{(by \ref{MS2b}, see \fullref{BasicMennicke}{Lam})}
 \\&= \Men{b g}{(f+ga)^2}{} \Men{q}{(f+ga)^2}{}
  && \text{(since $(f+ga)^2 \equiv 1 \mod q A$)}
 \\&= \Men{b}{(f+ga)^2}{} \Men{gq}{(f+ga)^2}{}
  && \text{(by \ref{MS2a})}
 \\&= \Men{b}{f+ga}{}^2
  && \begin{matrix}
 \text{(by \ref{MS2b} and because} \hfill
 \\ \text{\qquad $(f+ga)^2 \equiv 1 \mod gq A$)}
 . \end{matrix}
 \end{align*}
 \end{proof}

\subsection{Nonstandard analysis} \label{NonstandardSect}

 \begin{rem} \label{NonstandRem}
 Many of the results and proofs in \S\ref{3.}  use the theory of nonstandard
analysis, in the language and notation of \cite{StroyanLuxemburg}. This
enables us to express some of the arguments in a form that is less complicated
and more intuitive. In particular, it is usually possible to eliminate phrases
of the form ``for every ideal~$\Iq$, there exists an ideal~$\Iq'$," because the
nonstandard ideal~$\II$ \seeDefn{ODefn} can be used as~$\Iq'$ for any choice
of the ideal~$\Iq$ of~$A$. (Thus, $\II$ plays a role analogous to the
set of infinitesimal numbers in the nonstandard approach to Calculus.)

 As an aid to those who prefer classical proofs, Remark~\ref{NowNonstandRem}
provides classical reformulations of the nonstandard results. It is not
difficult to prove these versions, by using the nonstandard proofs as detailed
hints. Doing so yields a proof of Theorem~\ref{BddGen2} without reference to
nonstandard analysis.

The unpublished manuscript \cite{CKP} uses nonstandard models much more
extensively than we do here, in place of the Compactness Theorem
\pref{Cpctness}, for example \cf{2.1}. We have employed them only where they
have the most effect.
 \end{rem}

\begin{notation}[{(cf.\ \cite{StroyanLuxemburg})}] \ 
 \begin{itemize}
 \item For a given ring~$A$, we use $\*A$ to denote a (polysaturated)
nonstandard model of~$A$.
 \item If $X$ is an entity (such as an ideal, or other subset) that is
associated to~$A$, we use $\*X$ to denote the corresponding standard entity
of~$\*A$.
 \item For an element~$a$ of~$A$, we usually use~$a$ (instead of $\*a$) to
denote the corresponding element of~$\*A$.
 \end{itemize}
 \end{notation}
 
 Recall that the $*$-transform of a first-order sentence is obtained by replacing
 each constant symbol~$X$ with~$\*X$ \cite[Defn.~3.4.2, p.~27]{StroyanLuxemburg}.
 For example, the $*$-transform of $\forall a \in A, \exists b \in B, (a = b^2)$
 is $\forall a \in \*A, \exists b \in \*B, (a = b^2)$.
 
 \begin{Leibniz}[{}{\cite[(3.4.3), p.~28]{StroyanLuxemburg}}]
 A first-order sentence with all quantifiers bounded is true in~$A$
 if and only if its $*$-transform is true in~$\*A$.
 \end{Leibniz}

The following result of nonstandard analysis could be used in place of the
Compactness Theorem \pref{Cpctness} in our arguments.

\begin{lem}[{\cite[(2.1)]{CKP}}] \label{2.1}
 Suppose $G$ is a group and $\genset$ is a subset of~$G$. The following are
equivalent:
 \begin{enumerate}
 \item \label{2.1-1}
 $\genset$ boundedly generates $\gen{\genset}$;
 \item \label{2.1-2}
 $\gen{\*\genset} =\*{\gen{\genset}}$;
 \item \label{2.1-3}
 $\gen{\*\genset}$ is of finite index in $\*{\gen{\genset}}$;
 \item \label{2.1-4}
 there exists a $*$-finite subset $\Omega$ of $\*G$ with
$\*{\gen{\genset}} \subseteq \gen{\*\genset}\Omega$.
 \end{enumerate}
 \end{lem}

\begin{proof}
 ($\ref{2.1-1} \Rightarrow \ref{2.1-2}$)
 If $\genset$ boundedly generates $\gen{\genset}$, then there exists a
positive integer~$r$, such that $\gen{\genset} = \gennum{\genset}{r}$.
Then 
 $$\*{\gen{\genset}} = \*{\gennum{\genset}{r}} =
\gennum{\*\genset}{r} \subseteq \gen{\*\genset} .$$

($\ref{2.1-2} \Rightarrow \ref{2.1-3} \Rightarrow \ref{2.1-4}$) Obvious.

($\ref{2.1-4} \Rightarrow \ref{2.1-1}$) Let $\Omega' = \Omega \cap
\*{\gen{{\genset}}}$. Since $\Omega'$ is $*$-finite, there exists
$\omega \in \*\natural$, such that $\Omega' \subseteq
\*{\gennum{\genset}{\omega}}$. For any infinite $\tau \in \*\natural$,
we have $\gen{\*{\genset}} \subseteq \*{\gennum{\genset}{\tau}}$.
Therefore, letting $r = \omega + \tau$, we have ``There exists $r \in
\*\natural$, such that $\*{\gen{{\genset}}} =
\*{\gennum{\genset}{r}}$". By Leibniz' Principle, ``There exists $r \in \natural$,
such that $\gen{{\genset}} = \gennum{{\genset}}{r}$."
 \end{proof}

\subsection{Two results from number theory}

Our proofs rely on two nontrivial theorems of number theory. The first of
these is a version of Dirichlet's Theorem on primes in arithmetic
progressions. It is a basic ingredient in our arguments (cf.\
few generators property \pref{1.2-2}). The second theorem is used only to
establish the claim in the proof of Lemma~\ref{4.7}.

\begin{thm}[{\cite[(A.11), p.~84]{BMS}}] \label{DirichletDensity}
 Let
 \begin{itemize}
 \item $\ints$ be the ring of integers of an algebraic number field~$K$,
 and
 \item $\Norm \colon K \to \rational$ be the norm map.
 \end{itemize}
 For all nonzero $a,b \in \ints$, such that $a \ints + b \ints = \ints$, there exist
infinitely many $h  \in a + b \ints$, such that
 \begin{enumerate}
 \item $h \ints$ is a maximal ideal of~$\ints$,
 and
 \item $\Norm(h)$ is positive.
 \end{enumerate}
 \end{thm}

\begin{rem}
 The fact that $\Norm(h)$ can be assumed to be positive is not essential to any
of the arguments in this paper. However, it simplifies the proof of
Lemma~\fullref{4.3}{fg}, by eliminating the need to consider absolute values.
(Also, if $\Norm(h)$ were not assumed to be positive, then a factor of~$2$
would be lost, so $(16k)!$ would replace $(8k)!$ in the conclusion, but that
would have no impact on the main results.)
 \end{rem}

\begin{thm}[{\cite[p.~57]{Ostmann}}] \label{3primes}
 Let $r$ and~$m$ be any positive integers, such that $\gcd(r,m) = 1$. Then
there exists $M \in \integer$, such that if $t$ is an integer greater
than~$M$, and $t \equiv 3r \mod m$, then $t = p_1 + p_2 + p_3$, where each
$p_i$~is a rational prime that is congruent to~$r$ modulo~$m$.
 \end{thm}

We do not need the full strength of Theorem~\ref{3primes}, but only the
following consequence:

\begin{cor} \label{6primes}
 Let $r$ and~$m$ be any positive integers, with $\gcd(r,m) = 1$. If $t \in m
\integer$, then $t$ can be written in the form
 $$ t = p_1 + p_2 + p_3 - p_4 - p_5 - p_6 ,$$
 where each $p_i$~is a rational prime that is congruent to~$r$ modulo~$m$.
 \end{cor}

In fact, the arguments could be carried through with a weaker result that
uses more than~$6$ primes: if we assume only that every $t \in
m \integer$ can be written in the form
 $$ t = p_1 + p_2 + \cdots + p_c - p_{c+1} - p_{c+2} - \cdots - p_{2c} ,$$
 then the only difference would be that the constant $7k$ in the conclusion
of Lemma~\ref{4.7} would be replaced with $(2c+1)k$. This would have no
effect at all on the main results.

\section{First-order properties and bounded generation when $n \ge 3$}
\label{FirstOrderSect}

In \S\ref{GENsect} and \S\ref{EXPsect}, we define certain first-order
properties that any particular ring may or may not have. They are denoted
$\GEN(\ttt,\rr)$, and $\EXP(\ttt,\ell)$, for positive integers $\ttt$, $\rr$,
and~$\ell$.
 (In order to apply the Compactness Theorem \pref{Cpctness}, it is crucial
that, for fixed values of the parameters $\ttt$, $\rr$, and~$\ell$, 
these properties can be expressed by first-order sentences.) 
 We also show that the number rings~$\BS$ of interest to us satisfy these
properties for appropriate choices of the parameters (see 
\ref{4.4} and~\ref{4.5}).
 In \S\ref{1.}, we show that these properties (together with the stable range
condition $\SRquot$) imply that the order of the universal Mennicke group is
bounded \see{1.8}. Finally, in \S\ref{2.}, we establish that if $n \ge 3$, then
the elementary matrices boundedly generate a finite-index subgroup of
$\SL(n,\BS)$ \fullsee{BddGen>3}{LU}.

\begin{notation} \label{FirstOrderSectNotation}
 Throughout this section,
 \begin{itemize}
 \item $K$ is an algebraic number field,
 \item $k$ is the degree of~$K$ over~$\rational$,
 \item $\ints$ is the ring of integers of~$K$, 
 \item $B$~is an order in~$\ints$,
 \item $S$~is a subset of $B \smallsetminus \{0\}$ that is closed under multiplication,
 and
 \item $\Norm \colon K \to \rational$ is the norm map.
 \end{itemize}
 \end{notation}

\subsection{Few generators property $\GEN(\ttt,\rr)$} \label{GENsect}

We write down a simple first-order consequence of Dirichlet's Theorem
\pref{DirichletDensity} on primes in arithmetic progressions. It will be used
to bound the number of generators of the universal Mennicke group (see
Step~\ref{1.6} of the proof of Theorem~\ref{1.8}). In addition, the special case
$\GEN(2,1)$ also plays a key role in the proof of Proposition~\ref{MennInSL2}.

\begin{defn}[{\cite[(1.2)]{CKP}}] \label{1.2-2}
 For fixed positive integers $\ttt$ and~$\rr$, a commutative ring~$A$ is said
to satisfy $\GEN(\ttt,\rr)$ if and only if:
 for all $a,b \in A$, such that $aA + bA = A$, there exists $h \in a + b A$,
such that 
 $$ \text{$\displaystyle \frac{\U{ h A \vphantom{\big|}}}{\U{h A
\vphantom{\big|}}^\ttt}$ can be generated by $\rr$ or less elements.} $$
 (Recall that $\U{h A}$ denotes the group of units in $A/hA$.)
 \end{defn}

\begin{lem} \label{BSsurj}
 If $b \in B$ and $s \in S$, with $b \neq 0$, then
 \begin{enumerate}
 \item \label{BSsurj-B}
 $B \subseteq b \BS + s B$,
 \item \label{BSsurj-BS}
 $B + b \BS = \BS$,
 and
 \item \label{BSsurj-surj}
 the natural homomorphism from~$B$ to $\BS / b \BS$ is surjective.
 \end{enumerate}
 \end{lem}

\begin{proof}
 \pref{BSsurj-B} Because $B/bB$ is finite, and $\{s^n B\}$ is a decreasing
sequence of ideals, there exists $n \in \integer^+$, such that $bB + s^n B = b
B + s^{n+1} B$. Hence $s^n \in b B + s^{n+1} B$, so
 $$1 \in s^{-n} (b B + s^{n+1} B) = s^{-n} b B + s B \subseteq b \BS + s B .$$

 \pref{BSsurj-BS} For any $s_0 \in S$, we know, from~\pref{BSsurj-B}, that 
 $1 \in b \BS + s_0 B$. Therefore $1/s_0 \in b \BS + B$.

\pref{BSsurj-surj} This is immediate from \pref{BSsurj-BS}.
 \end{proof}

 \begin{prop}[{\cite[(4.1)]{CKP}}] \label{4.1} Let 
 \begin{itemize}
 \item  $a \in \BS$
 and $b \in B$,
 such that $b \neq 0$ and $a \BS + b \BS = \BS$,
 and
 \item $\gamma$ be any nonzero element of~$\ints$, such that $\gamma \ints \subset
B$.
 \end{itemize}
  Then:
 \begin{enumerate}
 \item \label{4.1-a0}
 there exists $a_0\equiv a \mod b \BS$, such that $a_0 B + b
\gamma^2 B = B$,
 and
 \item \label{4.1-a'}
 for any $a'$ in~$\ints$ with $a' \equiv a_0 \mod b\gamma^2 \ints$, 
 \begin{enumerate}
 \item \label{4.1-a'-homo}
 the natural homomorphism $\BS \to \localize{\ints}{S}/ a' \localize{\ints}{S}$ is
surjective, and has kernel~$a' \BS$,
 \item \label{4.1-a'-iso}
 $\BS/ a' \BS$ is isomorphic to a quotient of $\ints/a'\ints$,
 and
 \item \label{4.1-a'-inB}
  $a' \in B$.
 \end{enumerate}
 \end{enumerate}
 \end{prop}

\begin{proof}
 \pref{4.1-a0} From \fullref{BSsurj}{surj}, the natural homomorphism $B \to
\BS/b \BS$ is surjective, so we may choose $a_1 \in B$ with $a_1 \equiv  a
\mod b \BS$. Since $a$~is a unit in $\BS / b \BS$, then $a_1$~is a unit in
$B/(B \cap b \BS)$; thus, there exist $x \in B$ and $y \in B\cap b \BS$, such
that $a_1 x+y = 1$. Since $B/b \gamma^2 B$ is semi-local (indeed, it is
finite), there exists $a_0 \in a_1 + y B$, such that $a_0$ is a unit in $B/b
\gamma^2 B$. Then $a_0 \equiv a_1 \equiv a \mod b\BS$ and $a_0 B+b\gamma^2 B =
B$.

 \pref{4.1-a'-inB} We have $a' \in a_0 + b \gamma^2 \ints \subseteq B + b
\gamma B = B$.

\pref{4.1-a'-homo} Let $\varphi\colon \BS \to \localize{\ints}{S}/ a'
\localize{\ints}{S}$ be the natural homomorphism.
 Because $b \gamma^2 \ints \subseteq b
\gamma B$, we have 
 $$a' B + \gamma B \supseteq a' B + b \gamma B = a_0 B+b\gamma B \supseteq a_0
B+b\gamma^2 B = B .$$
 Therefore 
 $$a'\BS + \gamma \BS = \BS .$$
 Hence $a' \localize{\ints}{S} + \gamma \localize{\ints}{S} = \localize{\ints}{S}$ (which implies $a' \localize{\ints}{S} + \BS =
\localize{\ints}{S}$ --- that is, $\varphi$ is surjective). In other words, $a'$ is
relatively prime to~$\gamma$, so
 $$ a' \localize{\ints}{S} \cap \gamma \localize{\ints}{S} = a' \gamma \localize{\ints}{S} \subseteq a' \BS .$$
 The kernel of~$\varphi$ is
 \begin{align*}
 a'\localize{\ints}{S} \cap \BS
 &= a'\localize{\ints}{S} \cap (a'\BS+\gamma \BS)
 \\&= a'\BS + (a' \localize{\ints}{S} \cap \gamma \BS)
 \\&\subseteq  a' \BS + (a' \localize{\ints}{S} \cap \gamma
\localize{\ints}{S})
 \\&= a'\BS .
 \end{align*}

\pref{4.1-a'-iso} From \pref{4.1-a'-homo}, we see that
 $\BS/ a' \BS \iso \localize{\ints}{S} / a' \localize{\ints}{S} $.
On the other hand, the natural homomorphism $\ints \to \localize{\ints}{S} / a' \localize{\ints}{S}$ is
surjective \fullsee{BSsurj}{surj} and has $a' \ints$ in its kernel, so $\localize{\ints}{S} /
a' \localize{\ints}{S}$ is isomorphic to a quotient of $\ints/a'\ints$. 
The desired conclusion follows.
 \end{proof}

\begin{cor}[{(cf.\ \cite[(4.4)]{CKP})}] \label{4.4}
 $\BS$ satisfies $\GEN(\ttt,1)$, for every positive integer~$\ttt$.
 \end{cor}

\begin{proof}
 Proposition~\fullref{4.1}{a0} yields $a_0 \equiv a \mod b \BS$, such that
$a_0 B + b \gamma^2 B = B$. Then $a_0 \ints + b \gamma^2 \ints = \ints$, so
Dirichlet's Theorem \pref{DirichletDensity} yields $h \in a_0 + b
\gamma^2 \ints$, such that $h \ints$ is a maximal ideal. Therefore 
$\ints/h\ints$ is a finite field.

 From \fullref{4.1}{a'-iso}, we know that 
 $\BS/ h \BS$ is isomorphic to a quotient of $\ints/h\ints$; thus, $\BS/h \BS$
is either trivial or a finite field. In either case, the group of units is
cyclic, so the quotient $\U{h \BS}/\U{h \BS}^\ttt$ is also cyclic.
 \end{proof}

\subsection{Exponent property $\EXP(\ttt,\ell)$} \label{EXPsect}

We now introduce a rather technical property that is used to bound the
exponent of the universal Mennicke group (see Step~\ref{1.7} of the proof of
Theorem~\ref{1.8}).  Theorem~\ref{4.5} shows that this property holds in
number rings~$\BS$.

\begin{defn}[{\cite[(1.3)]{CKP}}] \label{1.3}
 Let $\ttt$ be a non-negative integer and let $\ell$ be a positive integer. A
commutative ring $A$ is said to satisfy $\EXP(\ttt,\ell)$ if and only if for
every $q$ in~$A$ with $q \neq 0$ and every $(a,b) \in W(q A)$,
there exists $a',c,d \in A$ and $u_i,f_i,g_i,b_i',d_i' \in A$ for $1 \le i \le
\ell$, such that
 \begin{enumerate} 
 \item \label{1.3-1}
 $a' \equiv a \mod bA$;
 \item \label{1.3-2}
 $\begin{bmatrix} a' & b \\ c & d \end{bmatrix}$  is in
$\SL(2,A; q A)$;
 \item \label{1.3-3}
 $\begin{bmatrix} a' & b_i' \\ c & d_i' \end{bmatrix}$
is in $\SL(2,A; q A)$ for $1 \le i \le \ell$;
 \item \label{1.3-4}
 $f_i I + g_i \begin{bmatrix} a' & b_i' \\ c & d_i' \end{bmatrix}$ is in
 $\SL(2,A; q A)$ for $1 \le i \le \ell$;
 \item \label{1.3-5}
 $(f_1+g_1 a')^2  (f_2+g_2 a')^2  \cdots (f_\ell+g_\ell a')^2 \equiv
(a')^{\ttt} \mod c A$;
 \item \label{1.3-6}
 $u_i$ is a unit in $A$ and $f_i+g_i a' \equiv u_i \mod  b_i' A$ for $1
\le i \le \ell$.
 \end{enumerate}
 \end{defn}

\begin{rem} \label{EXPeasy}
 Assume $\ttt$ is even, and let~$A$ be an arbitrary commutative ring.
 \begin{enumerate} \renewcommand{\theenumi}{\alph{enumi}}
 \item  \label{EXPeasy-almost}
 It is easy to satisfy all of the
conditions of Definition~\ref{1.3} \emph{except} the requirement that $u_1$~is a
unit: simply choose $f_1,g_1 \in A$, such that
 $$ \mat{a}{b}{c}{d}^{\ttt/2} = f_1 \Id{2} + g_1 \mat{a}{b}{c}{d} ,$$
 and let $f_i = 1$ and $g_i = 0$ for $i > 1$.
 \item \label{EXPeasy-a=0}
  If $a = 0$, then it is easy to satisfy all the conditions of Definition~\ref{1.3}:
  choose $f_i,g_i$ as in \pref{EXPeasy-almost}, and, because $b_i = b$ is
  a unit, we may let $u_i = 1$ for all~$i$.
 \item \label{EXPeasy-b=0}
 If $b = 0$, then it is easy to satisfy all the conditions of Definition~\ref{1.3}.
This is because $a$~must be a unit in this case, so we may let $u_1 = a^{\ttt/2}$
(and $u_i = 1$ for $i > 1$).
 \end{enumerate}
 \end{rem}

Recall that $k$ is the degree of~$K$ over~$\rational$
\see{FirstOrderSectNotation}.

\begin{lem}[{\cite[(4.3)]{CKP}}] \label{4.3}
 \ 
 \begin{enumerate}
 \item \label{4.3-homo}
 For any rational prime~$p$ and positive integer~$r$, let $\Norm_{p^r}$ be the
homomorphism from $\U{ p^r \ints}$ to $\U{ p^r \integer}$ induced
by the norm map $\Norm \colon K \to \rational$. If $p^r > 8k$, then the image
of $\Norm_{p^r}$ has more than~$2$ elements.
 \item \label{4.3-fg}
 If $f, g \in \BS$ with $f\BS+g\BS = \BS$, then, for any positive integer~$n$
and any nonzero $h \in \BS$, there exists $f' \equiv f \mod g \BS$, such that
 \begin{enumerate}
 \item \label{4.3-fg-e}
 $\gcd \bigl( \ee(f'\BS), n \bigr)$ is a divisor of\/ $(8k)!$, where
$\ee(f'\BS)$ is the exponent of $\U{f'\BS}$,
 and
 \item $f' \BS + h \BS = \BS$.
 \end{enumerate}
 \end{enumerate}
 \end{lem}

\begin{proof}
 \pref{4.3-homo}
 It is well known that $\U{p^r \integer}$ has a cyclic subgroup of order
$(p-1) p^{r-1}$ if $p$~is odd, or of order $p^{r-2}$ if $p = 2$. Thus, in any
case, $\U{p^r \integer}$ has a cyclic subgroup~$C$ of order $\ge p^r/4 > 2k$.
For $c \in \integer$ and, in particular, for $c \in C$, we have $\Norm(c) =
c^k$. Therefore
 $$ \# \Norm_{p^r} \bigl( \U{p^r \ints} \bigr)  \ge  \# \Norm_{p^r}(C) 
 \ge \frac{\# C}{\gcd \bigl( k, \# C \bigr)} > \frac{2k}{k} = 2 .$$

 \pref{4.3-fg} We may assume $h = n$, by replacing $n$ with $n \, |\Norm(h
s)|$, for some $s \in S$ with $h s \in B$.
 We consider two cases.

\setcounter{case}{0}

\begin{case} \label{4.3pf-BS=O}
 Assume $\BS = \ints$.
 \end{case}
 Choose $f_0 \equiv f \mod g \ints$, such that $f_0 \ints +  (8k)! \, n g \ints =
\ints$. 

 Let $P$ be the set of rational prime divisors of~$n$.
 We may assume (by replacing~$n$ with the product $\Norm(g) \, n$) that
$P$~contains every prime divisor of~$\Norm(g)$.
 For each $p$ in~$P$, let 
 $$ \text{$r(p)$ be the largest integer such that $p^{r(p)}$
divides~$(8k)!$.}$$ 
 From \pref{4.3-homo}, we know that the image of $\Norm_{p^{r(p)+1}}$ has more
than~$2$ elements. Therefore, $\Norm_{p^{r(p)+1}}(f_0)$ (or any other element
of the image) can be written as a product of two elements of the image, neither
of which is trivial.
 This implies that there exist $x(p),y(p) \in \ints$, such that 
 \begin{itemize}
 \item $x(p) \, y(p) \equiv f_0 \mod p^{r(p)+1} \ints$
 and
 \item neither $\Norm\bigl( x(p) \bigr)$ nor $\Norm\bigl( y(p) \bigr)$
is is congruent to $1$ modulo $p^{r(p)+1}$.
 \end{itemize}
 Now, by Dirichlet's Theorem~\pref{DirichletDensity} and the Chinese
Remainder Theorem, pick 
 \begin{itemize}
 \item $f_1 \in \ints$, such that 
 \begin{itemize}
 \item $f_1 \equiv x(p) \mod p^{r(p)+1} \ints$, for each $p$ in~$P$,
 \item $f_1 \ints$ is maximal,
 \item $n \notin f_1 \ints$,
 and
 \item $\Norm(f_1) > 0$;
 and
 \end{itemize}
 \item  $f_2 \in \ints$, such
that 
 \begin{itemize}
 \item $f_1 f_2 \equiv f_0 \mod \left(  \prod_{p \in P} p^{r(p)+1} \right) g \, \ints$,
 \item $f_2 \ints$ is maximal,
 \item $f_1 \ints + f_2 \ints = \ints$,
 and
 \item $\Norm(f_2) > 0$.
 \end{itemize}
 \end{itemize}
 Set $f' = f_1 f_2$, so $f' \equiv f_0 \equiv f \mod g \ints$. 

 The Chinese Remainder Theorem implies that $\U{f_1 f_2 \ints} \iso \U{f_1 \ints}
\times \U{f_1 \ints}$. Also, since $f_j \ints$ is maximal, for $j = 1,2$, we  know
that $\U{f_j \ints}$ is cyclic of order
 $$ \# \U{f_j \ints} = \# (\ints/ f_j \ints) - 1 = \# \Norm(f_j) - 1 =
\Norm(f_j) - 1 .$$
 Therefore
 $$ \ee(f' \ints) = \ee(f_1 f_2 \ints) = \lcm \bigl( \ee(f_1 \ints), \ee(f_1
\ints) \bigr)
 =  \lcm \bigl( \Norm(f_1) - 1, \Norm(f_2) - 1 \bigr) .$$
 For each $p$ in~$P$, we have
 $$ f_1 f_2 \equiv f_0 \equiv x(p) \, y(p) \equiv f_1 \, y(p) \mod p^{r(p)+1}
\ints, $$
 so $f_2 \equiv y(p) \mod p^{r(p)+1} \ints$. Thus, our selection of $x(p)$
and~$y(p)$ guarantees that 
 $$\text{$\gcd \bigl( \Norm(f_j) - 1 , n \bigr)$ is a divisor of~$(8k)!$ for
$j = 1, 2$.} $$ 
 Therefore 
 $$\gcd \bigl( \ee(f' \ints), n \bigr) = \lcm \Bigl( \gcd \bigl( \Norm(f_1) - 1
, n \bigr), \gcd \bigl( \Norm(f_2) - 1 , n \bigr) \Bigr) $$
 is a divisor of $(8k)!$.

\begin{case}
 The general case. 
 \end{case}
 We may assume $g \in B$, by replacing $g$ with $s g$, for some appropriate $s
\in S$. (Note that, since elements of~$S$ are
units in~$\BS$, we have $s g \BS = g \BS$.) Let $\gamma$ be a nonzero element
of~$\ints$, such that $\gamma \ints \subseteq B$. By \fullref{4.1}{a0}, there
exists $f_0 \in B$, such that $f_0 \equiv f \mod g \BS$ and $f_0 B + g n
\gamma^2 B = B$. From Case~\ref{4.3pf-BS=O}, we get $f' \equiv f_0 \mod g n
\gamma^2 \ints$, such that 
 $$ \text{$\gcd  \bigl( \ee(f' \ints), n \bigr)$ is a divisor of $(8k)!$.} $$
 From \fullref{4.1}{a'-iso}, we see that $\U{f' \BS}$ is isomorphic to a
quotient of $\U{f' \ints}$, so \pref{4.3-fg-e} holds. Also, we have $f' \equiv
f_0 \equiv f \mod g \BS$ and 
 \begin{align*}
 f' \BS + n \BS
 &\supseteq f' \BS + g n \gamma^2 \localize{\ints}{S} 
 = f_0 \BS + g n \gamma^2 \localize{\ints}{S}
 \\&\supseteq  f_0 \BS + g n
\gamma^2 \BS
 =  \BS .
 \end{align*}
 \end{proof}

\begin{thm}[{(cf.\ \cite[(4.5)]{CKP})}] \label{4.5}
 $\BS$ satisfies $\EXP \bigl(2(8k)!, 2 \bigr)$.
 \end{thm}

\begin{proof}
 Let 
 \begin{itemize}
 \item $q$ be any element of~$\BS$ with $q \neq 0$, 
 \item $(a,b)$ be an arbitrary element of $W(q \BS)$ with 
 $a \neq 0$ and $b \neq 0$ 
 (see~\ref{EXPeasy}(\ref{EXPeasy-a=0},\ref{EXPeasy-b=0})),
 \item $c,d \in \BS$, such that $\mat{a}{b}{c}{d} \in \SL(2,\BS; q \BS)$,
 \item $a' = a$, $b_1' = b$, $d_1' = d$, 
 \item $u_i = 1$ for $i = 1,2$,
 \item $b_0 = b/q \in \BS$,
 \item $\alpha_1$ be the exponent of~$a$ modulo $b\BS$,
 \item $b' \equiv b_0 \mod a \BS$, such that \fullsee{4.3}{fg}:
 \begin{itemize}
 \item the exponent $\alpha_2 = \ee( b' \BS)$ has the property that
$\gcd(\alpha_1,\alpha_2)$ is a divisor of $(8k)!$,
 and 
 \item $b' \BS + q \BS = \BS$,
 \end{itemize}
 \item $b_2' = b' q$,
 \item $d_2' \in \BS$, such that
 $\mat{a}{b_2'}{c}{d_2'} \in \SL(2,\BS; q \BS)$,
 that is,
 $$ d_2' = d + \frac{(b_2'-b)c}{a} = d + \frac{q (b'-b_0) c}{a} ,$$
 \item $t_1,t_2 \in \integer$, such that $\alpha_1 t_1 + \alpha_2 t_2 =
(8k)!$,
 and
 \item $f_i,g_i \in \BS$ (for $i = 1,2$) be defined by
 \begin{align*}
 \mat{a}{b}{c}{d}^{\alpha_1 t_1} &= f_1 \Id{2} + g_1 \mat{a}{b}{c}{d} , \\
 \mat{a}{b_2'}{c}{d_2'}^{\alpha_2 t_2} &= f_2 \Id{2} + g_2
\mat{a}{b_2'}{c}{d_2'} 
 . \end{align*}
 \end{itemize}
 Now, by multiplying matrices modulo $c \BS$, we see that 
 $$a^{2 (8k)!} = (a^{\alpha_1 t_1})^2 (a^{\alpha_2 t_2})^2  \equiv (f_1+g_1
a)^2 (f_2+g_2 a)^2 \mod c \BS .$$
 Similarly, 
 $$f_1+g_1 a \equiv a^{\alpha_1 t_1} \equiv 1 = u_1 \mod b \BS .$$
 Finally, because
 \begin{itemize}
 \item $f_2+g_2 a \equiv a^{\alpha_2 t_2} \mod b_2' \BS$ (by a similar
calculation),
 \item $a^{\alpha_2 t_2} \equiv 1^{t_2} = 1 \mod b' \BS$ (by definition
of~$\alpha_2$),
 \item $a^{\alpha_2 t_2} \equiv 1^{\alpha_2 t_2} = 1 \mod q \BS$ (since $(a,b)
\in W(q \BS)$),
 and
 \item $b_2' = b' q$ with $b'$ relatively prime to~$q$,
 \end{itemize}
 we conclude that 
 $$ f_2+g_2 a \equiv a^{\alpha_2 t_2} \equiv 1 = u_2 \mod b_2' \BS .$$
 \end{proof}

\begin{rem} \label{8k!->Const}
 The function~$2(8k)!$ in the conclusion of Theorem~\ref{4.5} is much larger than
necessary, but reducing the order of magnitude would not yield any
improvement in the main results --- all that matters is that the function
depends only on~$k$.
 However, it would be of interest to replace $2(8k)!$ with a function that is
bounded on an infinite subset of~$\natural$. For example, perhaps there is a
constant~$\ttt$ (independent of~$k$), such that $\BS$ satisfies $\EXP(\ttt,2)$
whenever $k$~is odd. If so, then the bound~$r$ in Theorem~\ref{EijBddGen} could
be chosen to depend only on~$n$, when $k$~is odd and $n \ge 3$.
 \end{rem}

\subsection{Bounding the order of the universal Mennicke group} \label{1.}
 The properties $\GEN(\ttt,\rr)$ and $\EXP(\ttt,\ell)$ were specifically
designed to be be what is needed in the proof of the following theorem.

\begin{thm}[{\cite[(1.8)]{CKP}}] \label{1.8}
 Let 
 \begin{itemize}
 \item $\ttt$, $\rr$, and~$\ell$ be positive integers,
 \item $A$ be an integral domain satisfying $\SRquot$, $\GEN(\ttt,\rr)$, and
$\EXP(\ttt,\ell)$,
 and
 \item $\Iq$ be an ideal in~$A$.
 \end{itemize}
 Then the universal Mennicke group $C(\Iq)$ is finite, and its order is bounded
by $\ttt^\rr$. 
 \end{thm}

\begin{proof}
 To bound the order of the abelian group~$C(\Iq)$, it suffices to bound both
the exponent and the number of generators needed. We assume $\Iq \neq 0$
(because the desired conclusion is obvious otherwise). 
 Note that, for any nonzero $q \in \Iq$, the natural homomorphism $C(q A) \to
C(\Iq)$ is surjective \fullsee{BasicMennicke}{smaller}, so we may assume $\Iq
= q A$ is principal.

\setcounter{step}{0}

\begin{step} 
 \label{1.7}
 {\rm (cf.\ \cite[(2.4)]{Liehl-Beschrankte})}
 The exponent of $C(\Iq)$ is a divisor of $\ttt$. {\rm(}I.e., if $z$ is
in $C(\Iq)$, then $z^{\ttt} = 1$.{\rm)} 
 \end{step}
 Let $\UMen{b}{a}{q A}$ be an arbitrary element of $C(q A)$. Because, by
assumption, $A$ satisfies the exponent property $\EXP(\ttt,\ell)$, there exist
$a',c,d \in A$ and $u_i,f_i,g_i,b_i',d_i' \in A$ (for $1 \le i \le \ell$)
satisfying the conditions of \pref{1.3}. Applying, in order,
 \fullref{1.3}{1}+\pref{MS1b},
 \fullref{AddlMen}{inverse},
 \pref{MS2b}+\fullref{1.3}{5},
 \fullref{AddlMen}{1.4},
 \fullref{AddlMen}{inverse},
 \fullref{AddlMen}{1.4},
 \fullref{1.3}{6}+\pref{MS1b},
 and
 \fullref{BasicMennicke}{unit},
 we have 
 \begin{align*}
 \UMen{b}{a}{q A}^{-\ttt}
 &= \UMen{b}{a'}{q A}^{-\ttt}
 = \UMen{c}{a'}{q A}^{\ttt}
 = \prod_i \UMen{c}{f_i + g_i a'}{q A}^2 
 \\&= \prod_i \UMen{c g_i}{f_i + g_i a'}{q A}^2 
 = \prod_i \UMen{b_i' g_i}{f_i + g_i a'}{q A}^{-2} 
 \\&= \prod_i \UMen{b_i'}{f_i + g_i a'}{q A}^{-2} 
 = \prod_i \UMen{b_i'}{u_i}{q A}^{-2} 
 = \prod_i 1^{-2} 
 = 1 
 . \end{align*}

\begin{step} 
 \label{1.6}
 $C(\Iq)$ can be generated by $\rr$ or less elements.
 \end{step}
 Because of Step~\ref{1.7}, it suffices to show, for each prime divisor~$p$
of~$\ttt$, that the rank of $C(\Iq)/C(\Iq)^p$ is $\le \rr$. 
 \begin{itemize}
 \item Let $\UMen{b_i}{a_i}{q A} \in C(q A)$, for $1 \le i \le \rr+1$.
 \item By repeated application of $\SRquot$, we can inductively construct a
sequence $a'_1,\ldots,a'_{\rr+1}$ of elements of~$A$, such that 
 \begin{itemize}
 \item $a_i' \equiv a_i \mod b_i A$, for $1 \le i \le \rr+1$,
 and
 \item $a_i' A + a_j' A = A$ for $1  \le i < j \le \rr+1$.
 \end{itemize}
 \item By the Chinese Remainder Theorem, choose $y \in A$ with 
 \begin{itemize}
 \item $y \equiv 1 \mod
q A$
 and
 \item $y \equiv b_i \mod a_i' A$ for $1 \le i \le \rr+1$.
 \end{itemize}
 \item Now $y A + a_1' a_2' \cdots a_{\rr+1}' q A = A$, so $\GEN(\ttt,\rr)$
implies that there exists 
 $$ \text{$h \equiv y \mod a_1' a_2' \cdots
a_{\rr+1}' q A$, such that $\U{h A}/\U{hA}^p$ has rank $\le \rr$.} $$
 \item Hence, there exists $\alpha \in A$ and integers $e_1, e_2, \ldots,
e_{\rr+1}$, with $e_i \not\equiv 0 \mod p$ for some~$i$, such that 
 $$ \prod_{i=1}^{\rr+1} (a'_i)^{e_i}  \equiv \alpha^p \mod h A .$$
 \item Since $h \equiv y \equiv 1 \mod q A$, we have $hA + q A = A$.
Therefore, we can choose $\beta \in A$ with 
 $$ \text{$\beta \equiv \alpha \mod h A$
 \quad and \quad
 $\beta \equiv 1 \mod q A$.} $$
  \end{itemize}
 We have $\beta^p \equiv \alpha^p \mod h A$ and $\beta^p \equiv 1 \mod q A$.
Hence
 $$\beta^p \equiv \prod_{i=1}^{\rr+1} (a'_i)^{e_i} \mod  hq A
.$$
 Now, since $h \equiv y \equiv b_i \mod a'_i A$, we have
 $$ \UMen{b_i}{a_i}{q A}
 = \UMen{b_i}{a_i'}{q A}
 = \UMen{b_i}{a_i'}{q A} \UMen{q}{a_i'}{q A}
 = \UMen{b_i q}{a_i'}{q A}
 = \UMen{hq}{a_i'}{q A} .$$
 Hence
 \begin{align*}
 \prod_{i=1}^{\rr+1} \UMen{b_i}{a_i}{q A}^{e_i}
 &= \prod_{i=1}^{\rr+1} \UMen{hq}{a_i'}{q A}^{e_i}
 = \UMen{hq}{\prod_{i=1}^{t+1}(a_i')^{e_i}}{q A}
 \\&= \UMen{hq}{\beta^p}{q A}
 = \UMen{hq}{\beta}{q A}^{p}
 \in C(q A)^p .
 \end{align*}
 Since  
 $\UMen{b_1}{a_1}{q A}, \ldots, \UMen{b_{\rr+1}}{a_{\rr+1}}{q A}$ are arbitrary
elements of~$C(q A)$, and some $e_i$~is nonzero modulo~$p$, we conclude that
the rank of $C(q A)/C(q A)^p$ is $\le \rr$.
 \end{proof}

\subsection{Bounded generation in $\SL(n,A)$ for $n \ge 3$} \label{2.}

The preceding results enable us to establish Theorem~\ref{EijBddGen} in the
case where $n \ge 3$ \fullsee{BddGen>3}{LU}.

\begin{thm} \label{SLnq/EnqFinite}
 Let 
 \begin{itemize}
 \item $n \ge 3$,
 \item $\ttt$, $\rr$, and~$\ell$ be positive integers,
 \item $A$ be an integral domain satisfying $\SRquot$, $\GEN(\ttt,\rr)$, and
$\EXP(\ttt,\ell)$,
 and
 \item $\Iq$ be an ideal in~$A$.
 \end{itemize}
 Then $\SL(n,A;\Iq)/\Enorm(n,A;\Iq)$ is finite, and its order is bounded by
$\ttt^\rr$. 
 \end{thm}

\begin{proof}
 By combining Thms.~\ref{BassThm} and~\ref{MennickeThm} (with $\mm = 2$ and $N
= \Enorm(n,A;\Iq)$), we see that $\SL(n,A;\Iq)/\Enorm(n,A;\Iq)$ is isomorphic
to a quotient of the universal Mennicke group $C(\Iq)$. 
 From Theorem~\ref{1.8}, we know that $\#C(\Iq) \le \ttt^\rr$, so the desired
conclusion is immediate.
 \end{proof}

Applying the Compactness Theorem \see{CpctnessCor} to this finiteness result
yields bounded generation. In the particular case of number rings, we obtain
the following conclusions.

\begin{cor}[{(cf.\ \cite[(2.4)]{CKP})}] \label{BddGen>3}
 Let
 \begin{itemize}
 \item $n$ be a positive integer $\ge 3$,
 \item $K$ be an algebraic number field,
 \item $k$ be the degree of~$K$ over~$\rational$,
 \item $B$ be an order in~$K$,
 \item $S$ be a multiplicative subset of~$B$,
 and
 \item $\Iq$ be an ideal in~$\BS$.
 \end{itemize}
 Then:
 \begin{enumerate}
 \item \label{BddGen>3-LU}
 $\LU(n,\BS)$ boundedly generates $\Elem(n,\BS)$,
 and
 \item \label{BddGen>3-LUnorm}
 $\LUnorm(n,\BS;\Iq)$ boundedly generates $\Enorm(n,\BS;\Iq)$.
 \end{enumerate}
 More precisely, there is a positive integer~$r$, depending only on~$k$
and~$n$, such that 
 $$ \text{$\gennum{\LU(n,\BS)}{r} = \Elem(n,\BS)$ and
$\gennum{\LUnorm(n,\BS;\Iq)}{r} = \Enorm(n,\BS;\Iq)$} .$$
 \end{cor}

A generalization of \fullref{BddGen>3}{LUnorm} that applies to all normal
subgroups, not merely the one subgroup $\Enorm(n,\BS;\Iq)$, can be found 
in \S\ref{NormalSect}.
It is proved by combining this result (and an analogous result for the case $n
= 2$) with the Sandwich Condition \pref{Sandwich}.

\section{Additional first-order properties of number rings}
\label{AddlPropSect}

We define two properties ($\UNIT(\rrr,\xx)$ and $\CONJ(\zz)$), and show they
are satisfied by number rings $\BS$ that have infinitely many units
\seeand{4.6}{4.7}.
 As in \S\ref{FirstOrderSect}, it is crucial that these properties can be
expressed by first-order sentences (for fixed values of the parameters $\rrr$,
$\xx$, and~$\zz$).

\begin{notation}
 Throughout this section,
 \begin{itemize}
 \item $K$ is an algebraic number field,
 \item $k$~is the degree of~$K$ over~$\rational$,
 \item $\ints$ is the ring of integers of~$K$, 
 \item $B$~is an order in~$\ints$,
 and
 \item $S$~is a multiplicative subset of~$B$.
 \end{itemize}
 \end{notation}

\subsection{Unit property $\UNIT(\rrr,\xx)$}

\begin{notation}
 If $u$ is a unit in a ring $A$, then
 $H(u)= \begin{bmatrix} u & 0 \\ 0 & u^{-1} \end{bmatrix}$.
 \end{notation}

\begin{defn}[{(cf.\ \cite[(3.1)]{CKP})}] \label{3.1}
 Let $\rrr$ and $\xx$ be positive integers. A commutative ring~$A$ satisfies
the \emph{unit property} $\UNIT(\rrr,\xx)$ if and only if:
 \begin{enumerate}
 \item \label{3.1-1}
 for each nonzero $q \in A$, there exists a unit~$u$ in~$A$, such that $u
\equiv 1 \mod q A$ and $u^4 \neq 1$;
 and 
 \item \label{3.1-3}
 there exists a unit $u_0$ in~$A$ with $u_0^2 \neq 1$, such that whenever
 \begin{itemize}
 \item $\Iq$ is an ideal in~$A$ with $\le \rrr$ generators,
 and
 \item $T \in
\SL(2,A;\Iq)$,
 \end{itemize}
  there exist $E_1,E_2,\ldots,E_\xx \in \LU(2,\Iq)$, such that
 $$H(u_0)^{-1} \, T \, H(u_0) = E_1 T E_2 E_3  \cdots E_\xx .$$
 \end{enumerate}
 \end{defn}

\begin{lem}[{\cite[(4.6)]{CKP}}] \label{4.6}
 If $\BS$ has infinitely many units, then $\BS$ satisfies the unit property
$\UNIT(\rrr,5)$, for any~$\rrr$.
 \end{lem}

\begin{proof}
 Let $q$ be any nonzero element of~$\BS$. Since, by assumption, $\BS$ has
infinitely many units, there is a unit $u$ in~$\BS$ that is not a root of
unity. Some power of~$u$ satisfies the requirements of \fullref{3.1}{1}.

  Let $u_0 = u^{(8k)!}$, 
 $\Iq$ be any ideal of~$\BS$, and 
 $T = \mat{a}{b}{c}{d}$ be any element of $\SL(2,A;\Iq)$. 
 We may assume $\Iq \neq 0$, for otherwise $T = \Id{2}$, so the
 conclusion of \fullref{3.1}{3} is trivially true.
 \begin{itemize}
 \item By \fullref{4.3}{fg-e}, there exists $a' \equiv a \mod b^2 \BS$,
such that 
 $$ \text{$\gcd \bigl( \ee(a \BS), \ee(a'\BS) \bigr)$ is a divisor of
$(8k)!$,} $$
 where $\ee(a \BS)$ denotes the exponent of $\U{a \BS}$. 
 \item Choose $z \in \BS$, such that $a' = a + z b^2$.
 \item Choose $t,t' \in \BS$, such that $t \ee(a \BS) + t' \ee(a' \BS) =
(8k)!$. 
 \item Let $u^{2t\ee(a \BS)}=ax+1$ and $u^{-2t \ee(a \BS)} =a y+1$. 
 \end{itemize}
 We have
 \begin{align*}
 H(u^{-t\ee(a \BS)}) T \, H(u^{t\ee(a \BS)})
 &= \mat{a}{u^{-2t\ee(a \BS)}b}{u^{2t\ee(a \BS)}c}{d}
 \\&= \mat{a}{(a y + 1)b}{(ax + 1)c}{d}
 \\&= \mat{1}{0}{x c}{1}
 \mat{a}{b}{c}{d}
 \mat{1}{y b}{0}{1}
 \\&=  E_{2,1}(*) \, T \, E_{1,2}(*)
 . \end{align*}
 Letting 
 $T' = T \mat{1}{0}{z b}{1} = \mat{a'}{b}{*}{d} $,
 the same calculation shows that 
 $$ H(u^{-t'\ee(a' \BS)}) T' \, H(u^{t'\ee(a' \BS)})
 =  E_{2,1}(*) \, T' \, E_{1,2}(*)
 . $$
 Therefore
 \begin{align*}
 H(u_0)^{-1} T \,  H(u_0)
 \hskip -0.75in &
 \\&= H(u^{-t' \ee(a' \BS)}) \,
 E_{2,1}(*) \,
 T \,
 E_{1,2}(*) \,
 H(u^{t' \ee(a' \BS)})
 \\&= H(u^{-t' \ee(a' \BS)}) \,
 E_{2,1}(*) \,
 T' \,
 E_{2,1}(*) \,
 E_{1,2}(*) \,
 H(u^{t' \ee(a' \BS)})
 \\&= 
 E_{2,1}(*) \,
 T' \, E_{1,2}(*) \,
 E_{2,1}(*) \,
 E_{1,2}(*) 
 \\&= 
 E_{2,1}(*) \,
 T \,
 E_{2,1}(*) \,
 E_{1,2}(*) \,
 E_{2,1}(*) \,
 E_{1,2}(*) 
 , \end{align*}
 with each of the elementary matrices in $\LU(2,\Iq)$.
Hence, \fullref{3.1}{3} is satisfied with $\xx = 5$.
 \end{proof}

\subsection{Conjugation property $\CONJ(\zz)$}

The following property will be used to control the image of the other
elementary matrices under conjugation by $E_{1,2}$ \see{4.7easy}.

\begin{defn} \label{3.2?}
 Let $\zz$ be a positive integer, and let $A$ be a commutative ring.
 \begin{itemize}
 \item For ideal~$\Iq$ of~$A$, let
 $$ \lowset_{\Iq} = \bigset{ y \in \Iq }{
 \begin{matrix} \text{there exists $z \equiv \pm1 \mod \Iq$, and}
 \\ \text{units $u_1,u_2$ in~$A$, such that}
 \\ \text{$1 + y z u_1^2 = u_2^2$} 
 \end{matrix} }.$$
 \item The ring~$A$ is said to satisfy $\CONJ(\zz)$ if, for every
nonzero $q \in A$, there is a nonzero $q' \in A$, such that every
element of~$q' A$ is a sum of $\le \zz$ elements of~$\lowset_{q A}$.
 \end{itemize}
 \end{defn}

 Most of the proof of the following theorem appears in
\cite[p.~327]{Vaserstein-SL2} and \cite[pp.~519--521]{Liehl-Beschrankte}, but
\cite{CKP} modified the argument to avoid Liehl's assumption that the
prime~$p$ splits completely in~$K$. This eliminates the need to place
restrictions on~$K$ (as in \cite{Liehl-Beschrankte}).

\begin{thm}[{(cf.\ \cite[(4.7)]{CKP})}] \label{4.7}
 If $\BS$ has infinitely many units, then $\BS$ satisfies $\CONJ(7k)$.
 \end{thm}

\begin{proof}
 For convenience, let us use $[j] \lowset_{q A}$ to denote the set of elements
of~$\BS$ that are a sum of $j$~elements of~$\lowset_{q A}$.

\begin{claim*}
 It suffices to find
 \begin{itemize}
 \item nonzero $q', q'' \in A$,
 \item positive integers $r,m$, with $\gcd(r,m) = 1$,
 and
 \item a finite subset $D$ of~$B$, such that\/ $\#D \le k + 1$, and the
$\integer$-span of~$D$ contains $q' B$,
 \end{itemize}
such that $p d q'' \in \lowset_{q A}$, for
 \begin{itemize}
 \item every rational prime $p$ that is congruent to~$r$ modulo~$m$,
 and
 \item every $d \in D$.
 \end{itemize}
  \end{claim*}
 \noindent 
 We show that if there exist such $q'$, $q''$, $r$, $m$, and~$D$, then the
principal ideal $q' q'' m \BS$ is contained in $[7k] \lowset_{q A}$. To
this end, let $b$ be any nonzero element of~$B$ and $s \in S$. By
assumption on~$D$, we may write 
 $q' b = \sum_{i=1}^{k+1} b_i d_i$ with $b_i \in \integer$ and $d_i \in D$.
For each i, 
 \begin{itemize}
 \item $b_i m$ is a signed sum of $6$~rational primes that are congruent
to~$r$ modulo~$m$ \see{6primes},
 \item $p d_i q'' \in \lowset{qA}$, for each of these primes~$p$,
 and
 \item $-\lowset_{q A} = \lowset_{q A}$ (because $z$ can be replaced by~$-z$),
 \end{itemize} so $b_i m d_i q''\in [6] \lowset_{q A}$.  Hence
 $$q'bmq''
 = \sum_{i=1}^{k+1} b_i m d_i q''
 \in [6(k+1)] \lowset_{q A} \subseteq  [7k] \lowset_{q A} .$$
 Since it is clear from the definition that $\lowset_{q A}$ is closed under
multiplication by~$s^{-2}$ (this is the reason for including the
unit~$u_1$), we see that $q'mq''(b/s^2) \in [7k] \lowset_{q A}$. Since
$b/s^2$ is an arbitrary element of~$\BS$, we conclude that $q'mq'' \BS
\subseteq [7k] \lowset_{q A}$, as desired.

This completes the proof of the claim. \qed

\medbreak

We now find $q'$, $q''$, $r$, $m$, and~$D$ as described in the Claim. We begin
by establishing notation.
 \begin{itemize}
 \item  To prove the result for a particular value of~$q$, it suffices to prove
it for some non-zero multiple of~$q$. Therefore, we may assume $q$ is a
rational integer, such that
 \begin{itemize}
 \item the exponent~$e$ of $\U{q \integer}$ is divisible by~$k!$,
 \item $q \ints \subseteq B$, 
 and
 \item  the discriminant of~$K$ divides~$q$.
 \end{itemize}
 \item Furthermore, we may assume there exists $t \in \integer$, such that
 $t^e - 1 \equiv a q \mod q^2 \integer$, with $\gcd(a,q) = 1$.
 (To achieve, this, let $q =  p_1^{e_1} p_2^{e_2} \cdots p_l^{e_l}$ be the
prime factorization. By carefully enlarging~$q$, we may assume, for $i \neq
j$, that $p_i^{e_i}$ does not divide $\phi(p_j^{e_j})$. Choose $t$ so that, for
each~$i$, its image in $\U{p_i^{e_i+2} \integer}$ is an element of maximal
order.)
 \item Let $D = \bigset{ t^{e-1}(a + q d_0) }{ d_0 \in D_0 \cup \{0\} }$,
where $D_0$ is some basis for~$B$ as a $\integer$-module. (Note that $\#D \le
k + 1$, and the $\integer$-span of~$D$ contains $t^{e-1} q B$.)
 \item Let $b = t^{e-1} a \prod_{d_0 \in D_0} (a + q d_0)$. (Note that every
element of~$D$ is a divisor of~$b$, and $b$~is relatively prime to~$q$.)
 \item Because $\BS$ has infinitely many units, there is some unit~$u$ that is
not a root of unity. Multiplying by an element of~$S$, we may assume $u \in
B$. Furthermore, by replacing~$u$ with an appropriate power, we may assume
that $u \equiv 1 \mod q^2 b B$ (and $u^2 \neq 1$). 
 \item Let $y = u^2-1 \in q^2 b B$, so $\Norm(y) \in (q^2 b B) \cap \integer$.
 \item Let $q'' = q y (1 + y)^{-1} = q y / u^2 \in \BS$.
 \item Write $\Norm(y) = n_0 n_1$, where $\gcd(q,n_0) = 1$, and any
(rational) prime dividing $n_1$ divides~$q$.
 \item Let $r$ be a rational integer with 
 \begin{itemize}
 \item $r \equiv t \mod q^2$
 and
 \item $r \equiv 1 \mod n_0$.
 \end{itemize}
 \end{itemize}

Let $d$ be any element of~$D$, and let $p$ be any rational prime that is
congruent to~$r$, modulo $\Norm(y)$.  We will show that $p d q'' \in
\lowset_{q A}$, which, by the Claim, completes the proof.
 \begin{itemize}
 \item  Since $y^{p^e - 1} \equiv 0 \mod dq^2 \BS$, and $p\equiv r \mod q^2 d
\BS$, and $td = t^e(a + q d_0) \equiv a \mod q \BS$, we have
 \begin{itemize}
 \item $y^{p^e - 1} + p^e - pdq - 1 \equiv 0 + r^e - 0 - 1 \equiv 0 + 1^e - 0 -
1 = 0 \mod d \BS$,
 and
 \item $y^{p^e - 1} +  p^e - pdq - 1 \equiv 0 + (t^e - 1) - tdq \equiv 0 + aq
- aq = 0 \mod q^2 \BS$,
 \end{itemize}
 so $y^{p^e - 1} + p^e - pdq - 1 \equiv 0 \mod dq^2 \BS$. Therefore
 $$ y^{p^e} + p^e y - pdqy - y $$
 is divisible by $dq^2y$.
 Since $k!$ divides~$e$ (and because $p$, being relatively prime to~$q$, does
not divide the discriminant of~$K$), we know that $y^{p^e} \equiv y \mod p
\BS$, so the displayed expression is also divisible by~$p$.
 \item Therefore $y^{p^e} + p^e y \equiv pdqy + y \mod pdq^2y \BS $.
 \item Hence $(1+y)^{p^e}
 \equiv 1 + p^e y + y^{p^e} 
 \equiv 1 + y + p d q y \mod pdq^2 y \BS$.
 \item Thus (recalling that $q'' = q y(1+y)^{-1}$), we have
 $$(1+y)^{p^e - 1}
 \equiv 1 + p d q'' \mod pdqq'' \BS ,$$
 so we may write
 $(1+y)^{p^e - 1} = 1 + pdq'' z$ with $z \equiv 1 \mod q \BS$.
 \item  Since $(1+y)^{p^e - 1} = (u^{p^e-1})^2$ is the square of a unit, we
conclude that $p d q'' \in \lowset_{q A}$ (taking $u_1 = 1$).
 \end{itemize}
 By the Claim, this completes the proof.
 \end{proof}

\section{Bounded generation in $\SL(2,A)$} \label{3.}

In this section, we establish Theorem~\ref{EijBddGen} in the case where $n = 2$
\fullsee{BddGen2}{LU}. This complements Theorem~\fullref{BddGen>3}{LU}, which
dealt with the case where $n \ge 3$.

\begin{rem}[(nonstandard analysis)] \label{NowNonstandRem}
 In this section, we frequently use the theory of nonstandard analysis
\cfSect{NonstandardSect}. As an aid to the reader who wishes to construct a
classical proof, we point out that:
 \begin{itemize}
 \item Corollary~\ref{3.6} is simply a restatement of Lemma~\ref{3.3} in nonstandard
terms.
 \item  Lemma~\ref{3.5} is a technical result that should be omitted from a
classical presentation of this material.
 \item  Proposition~\ref{3.7} asserts the existence of an ideal~$\Iq'$ of~$A$, such
that
 $$[\SL(2,K) ,\SL(2,A;\Iq') ] \subseteq \Elem(2,\Iq) .$$
 \item Lem~\ref{3.11} states, for any nonzero $y \in A$, that there is  a
nonzero ideal~$\Iq'$ of~$A$, such that $\Iq' \subseteq \Iq$, and
 $ \Men{by^2}{a}{\Iq} = \Men{b}{a}{\Iq} $
 for all $(a,b)\in W(\Iq')$.
 \end{itemize}
 \end{rem}

\begin{notation} \label{ODefn}
 Let
  $$\displaystyle \II
 = \mkern -10mu \bigcap_{\begin{matrix} 
 \text{$\Iq$ is an} \\ \text{ideal in~$A$} \end{matrix} } \mkern -20mu \*\Iq
\mkern 10mu
 = \bigcap_{\textstyle q \in A} \mkern -5mu q \, \*A .$$
 This is an (external) ideal of~$\*A$.
 \end{notation}

\subsection{Preliminaries}

\begin{defn}[{(cf.\ \cite{Vaserstein-SL2})}] \label{VasDefn}
 If $\Iq$ is an ideal in a commutative ring~$A$, then
 $$\Vas(2,A;\Iq) = \bigset{
 \begin{bmatrix} a & b \\ c & d \end{bmatrix} \in \SL(2,A)
 }{
 \begin{matrix}
 a,d \equiv 1 \mod \Iq^2 \\
 b,c \equiv 0 \mod \Iq
 \end{matrix}
 } .$$
 This is a subgroup of $\SL(2,A;\Iq)$ that contains $\Elem(2,\Iq)$ and
$\SL(2,A; \Iq^2)$.
 \end{defn}

\begin{lem} \label{3.5}
 The ideal $\II$ is nonzero, and we have $\II^2 = \II$.
 Therefore, 
 $$\Vas(2,A;\II) = \SL(2,A;\II) .$$
 \end{lem}

\begin{proof}
 If $F$ is any finite  set of nonzero elements of~$A$, then (because $A$ is an
integral domain) there is some nonzero $y \in A$, such that $y$~is a multiple
of every element of~$F$. Since $\*A$ is assumed to be polysaturated, this
implies that there is some nonzero $z \in \*A$, such that $z$~is a multiple of
every element of~$A$. Then $z \in \II$, so $\II \neq \{0\}$.

 Now, for any fixed element~$z$ of~$\II$, consider the internal binary
relation $R_z \subseteq \*A \times \*A$ given by 
 $$R_z = \{\, (x,y) \mid \text{$y \in \*A x$ and $z \in \*A y^2$} \,\}. $$
 Since $z \in \II$, it is easy to see that if $F$ is any finite set of nonzero
elements of~$A$, then there is a nonzero element~$y$ of~$A$, such that $(x,y)
\in R_z$ for all $x \in F$. By polysaturation, there is some nonzero $y_0 \in
\*A$, such that $(x,y_0) \in R_z$ for every nonzero~$x \in A$. Therefore $y_0
\in \II$ and $z \in \*A y_0^2 \subseteq \II^2$. Since $z$ is an arbitrary element
of~$\II$, we conclude that $\II^2 = \II$.
 \end{proof}

\begin{lem}[{(Vaserstein \cite{Vaserstein-SL2})}] \label{VasId}
  Let $A$ be a commutative ring and $u$~be a unit in~$A$.
 \begin{enumerate} \renewcommand{\theenumi}{A\arabic{enumi}}
 \item \label{VasId-A1}
 Suppose $u \equiv 1 \mod q^2A$ for some $q$ in~$A$. Then $u = 1 + xy$ with
$x,y$ in~$qA$, and we have
 $$ H(u) = \mat{1}{x}{0}{1} \mat{1}{0}{y}{1} \mat{1}{-u^{-1} x}{0}{1}
\mat{1}{0}{-uy}{1} .$$

\item \label{VasId-A2}
 $\mat{1}{1}{0}{1} H(u) \mat{1}{1-u^{-2}}{0}{1}
\mat{1}{-1}{0}{1} = H(u)$.

\item \label{VasId-A3}
 For $x$ in~$A$, we have
 $H(u)^{-1} \mat{1}{0}{x}{1} H(u) = \mat{1}{0}{x u^2}{1}$.

\item \label{VasId-A4}
 For $y,z \in A$, set
 $$ M(y,z) = 
 \mat{1}{z-1}{0}{1}
 \mat{1}{1}{0}{1}
 \mat{1}{0}{y}{1}
 \mat{1}{-1}{0}{1}
 .$$
 Suppose $u^2-1$ is in $(1 + y z)A$, and $(1 + yz)w = u^2-1$.
Set $c = w(1 - z + y z)$. Then
 $$M(y,z)H(u)^{-1} \mat{1}{c}{0}{1} M(y,z)^{-1} H(u) = \mat{1}{0}{-w y}{1} .$$

 \end{enumerate}
 \end{lem}

An argument similar to the proof of Corollary~\ref{VasLem1Normal} establishes the
following result.

\begin{lem}[{(Vaserstein's Lemma 1, cf.\ \cite[Lem.~1]{Vaserstein-SL2})}]
\label{VasersteinLemma}
 Let
 \begin{itemize}
 \item $A$ be a commutative ring,
 and
 \item $\Iq$ and~$\Iq'$ be nonzero ideals of~$A$, such that $\Iq' \subseteq
\Iq^2$.
 \end{itemize}
 If $A/\Iq'$ satisfies $\SR_1$, then $\Vas(2,A;\Iq) = \SL(2,A;\Iq')\Elem(2,
\Iq)$.
 \end{lem}

\begin{proof}
 By modding out $\Iq'$, we may assume that $A$ satisfies $\SR_1$, and we wish
to show that $\Elem(2, \Iq) = \Vas(2,A;\Iq)$. It suffices to
show that if $a,b \in A$, with
 \begin{enumerate}
 \item $a \equiv 1 \mod \Iq^2$,
 \item $b \in \Iq$,
 and
 \item $a A + b A = A$,
 \end{enumerate}
  then there exists $E \in \Elem(2, \Iq)$, such that $(a,b) E = (1,0)$.

Since $a \equiv 1 \mod \Iq$, and $a  A + b A = A$, we know that $aA + b\Iq =
A$; so $\SR_1$ implies that there exists $q \in \Iq$, such that $a + b q$ is
a unit. Thus, by replacing $(a,b)$ with $(a,b) E_{2,1}(q)$, we may assume
 $$ \text{$a$~is a unit.} $$
 Then, by replacing $(a,b)$ with $(a,b) E_{1,2}(-a^{-1} b)$, we may assume
 $$ b = 0 .$$
 Write 
 $$1 - a = x_1 y_1 + \cdots + x_r y_r$$
 with $x, y \in \Iq$, and $r$~minimal.
 The remainder of the proof is by induction on~$r$.

\emph{Base case. Assume $a = 1 + x y$ with $x, y \in \Iq$.}
  Applying $E_{1,2}\bigl( a^{-1} x  \bigr)$, $E_{1,2}( -y)$, and $E_{1,2}(-x)$
sequentially, we have
 $$(a,0) 
 \rightarrow (a, x) 
 \rightarrow (1,x) 
 \rightarrow (1,0) .$$

\emph{Induction step.} Let $\Iq'' = x_1 y_1 A + \cdots + x_{r-1} y_{r-1} A$.
 Now
 $$1 - a = x_1 y_1 + \cdots + x_r y_r \equiv x_r y_r \mod \Iq'' ,$$
 so, by applying the base case to the ring $A/\Iq''$, we know there is some $E
\in \Elem(2, \Iq)$, such that $(a,0)E \equiv (1,0) \mod \Iq''$. We may
also assume, by the argument above, that $(a,0)E = (u,0)$, for some
unit~$u$ (because the transformations will not change
the congruence class of $(a,0)E$ modulo~$\Iq''$). By the induction hypothesis,
then there exists $E' \in \Elem(2, \Iq)$, such that
 $(a,0) E E' = (1,0)$.
 \end{proof}

\subsection{A sufficient condition for a Mennicke symbol}

Because Mennicke's Theorem \fullref{MennickeThm}{MS2a} does not apply when $n =
2$, we prove the following result that yields a Mennicke symbol.

\begin{prop} \label{MennInSL2}
 Suppose
 \begin{itemize}
 \item $A$ is an integral domain,
 \item $\Iq$ is an ideal of~$A$,
 and
 \item $N$ is a normal subgroup of\/ $\SL(2,A;\Iq)$,
 \end{itemize}
 such that
 \begin{itemize}
 \item $A$ satisfies $\SRquot$ and $\GEN(2,1)$, 
 \item $\bigl[ \Elem(2,A), \SL(2,A;\Iq) \bigr] \subseteq N$,
 \item $\Enorm(2,A;\Iq) \subseteq N$,
 \item $C = \SL(2,A;\Iq)/N$,
 \item $\Men{\ }{\ }{} \colon W(\Iq) \to C$ is defined by $\Men{b}{a}{} =
\mat{a}{b}{*}{*} N$,
 and
 \item $\*{\Men{b y^2}{a}{}} = \*{\Men{b}{a}{}}$ for all $(a,b) \in W(\II)$
and all nonzero $y \in A$.
 \end{itemize}
 Then $\Men{\ }{\ }{}$ is a well-defined Mennicke symbol.
 \end{prop}

\begin{rem} \label{MS1OK}
 From Theorem~\ref{MennickeThm}, we know that $\Men{\ }{\ }{}$ is well defined,
and satisfies \pref{MS1a} and \pref{MS1b}. The problem is to establish
\pref{MS2a}.
 \end{rem}

We begin with a useful calculation.

\begin{lem}[{(cf.\ \cite[Case~1, p.~331]{Vaserstein-SL2})}] \label{MenMultInv}
 Let $(a_1,b),(a_2,b) \in W(\Iq)$, and suppose 
 $\mat{a_2}{b}{c}{d} \in \SL(2,A;\Iq)$. Then
 $$ \Men{b}{a_1a_2}{} \Men{b}{a_2}{}^{-1} 
 = \Men{a_2 b (1-a_1)}{1+a_2 d(a_1-1)}{} .$$
 \end{lem}

\begin{proof}
 Because $b c = a_2 d - 1$, we have
 \begin{align*}
 \mat{a_1 a_2}{b}{*}{*} \mat{a_2}{b}{c}{d}^{-1}
 &= \mat{a_1 a_2}{b}{*}{*} \mat{d}{-b}{-c}{a_2}
 \\&= \mat{a_1 a_2 d - b c}{-a_1 a_2 b + b a_2}{*}{*}
 \\&= \mat{1 + a_2 d (a_1 - 1)}{a_2 b(1 - a_1)}{*}{*}
 .
 \end{align*}
 \end{proof}

Let us show that it suffices to consider principal ideals.

\begin{lem}[{\cite[(3.13)]{CKP}}] \label{MenSL2CanPrinc}
  Suppose $A$, $\Iq$, $N$, $C$, and $\Men{\ }{\ }{}$ are as in the statement of
Proposition~\ref{MennInSL2}.
 If the restriction of $\Men{\ }{\ }{}$ to $W(q A)$ is a
Mennicke symbol, for every nonzero $q \in \Iq$, then $\Men{\ }{\ }{}$ is a
Mennicke symbol.
 \end{lem}

\begin{proof}
  By \pref{MS1OK} and Lam's Theorem~\fullref{MS2b->MS2a}{full}, we need
only establish \pref{MS2b}. Given 
 $$(a_1,b), (a_2,b) \in W(\Iq) ,$$
 we know, by assumption, that
 \begin{equation} \label{MenSL2CanPrincPf-restrict}
 \text{the restriction of $\Men{\ }{\ }{}$ to $W(q A)$ is a
Mennicke symbol}
 \end{equation}
 we have
 \begin{align*} \Men{b}{a_1 a_2}{} \Men{b}{a_2}{}^{-1}
 &= \Men{a_2 b (1-a_1)}{1 + a_2 d (a_1 - 1)}{} 
 &&\text{\pref{MenMultInv}}
 \\&= \Men{a_2 b (1-a_1)}{1 + a_2 d (a_1 - 1)}{}
 \Men{1-a_1}{1+a_2 d (a_1 - 1)}{}
 &&\text{(2nd is trivial)}
 \\&= 
 \Men{a_2 b (1-a_1)^2}{1 + a_2 d (a_1 - 1)}{}
 && \pref{MenSL2CanPrincPf-restrict}
 \\&= \Men{a_2 (1-a_1)}{1 + a_2 d (a_1 - 1)}{}
  \Men{b (1-a_1)}{1 + a_2 d (a_1 - 1)}{}
 && \pref{MenSL2CanPrincPf-restrict}
 \\&= 1 \cdot \Men{b (1-a_1)}{1 + (1 + bc) (a_1 - 1)}{}
 && \pref{MS1b}
 \\&= \Men{b (1-a_1)}{a_1 + bc (a_1 - 1)}{}
 \\&= \Men{b (1-a_1)}{a_1}{}
 && \pref{MS1b}
 \\&= \Men{b}{a_1}{}
 && \pref{MS1a}
 . \end{align*}
 \end{proof}

\begin{proof}[\bf Proof of Proposition~\ref{MennInSL2}]
 By Lemma~\ref{MenSL2CanPrinc}, we may assume $\Iq = q A$ is principal. Also,
by \pref{MS2bTriv->MS2a}, it suffices to show that if $(a_1,b q), (a_2,b q)
\in W(\Iq)$, and either $\Men{bq}{a_1}{} = 1$ or $\Men{bq}{a_2}{} = 1$, then
 \begin{equation} \label{PrincMennickeSL2PF-MS2b}
 \Men{b q}{a_1}{} \Men{b q}{a_2}{} = \Men{bq}{a_1 a_2}{}
 . \end{equation}
 Note that the elements $\Men{bq}{a_1}{}$ and $\Men{bq}{a_2}{}$ commute with
each other (because one of them is trivial). Thus, there is no harm in
interchanging $a_1$ with~$a_2$ if it is convenient. (That is why we do not
assume it is $\Men{bq}{a_2}{}$ that is trivial; it is better to allow
ourselves some flexibility.)

\setcounter{case}{0}

\begin{case} \label{3.15}
 Assume that either $a_1$ or~$a_2$ is a square modulo~$b q A$.
 \end{case}
 Because there is no harm in interchanging $a_1$ with~$a_2$, we may assume it
is~$a_2$ that is a square modulo~$b q A$.

Applying \pref{MS1a} and \pref{MS1b} allows us to make some
simplifying assumptions:
 \begin{itemize}
 \item By adding a multiple of $bq$ to~$a_2$, we may assume $a_2 = y^2$, for
some $y$ in~$A$.
 \item By adding a multiple of~$a_1 a_2$ to~$b$, we may assume
$b \*A + \II = \*A$ (because $\*A/\II$ satisfies $\SR_1$).
 \item By adding a multiple of~$bq$ to~$a_1$, we may assume $a_1 \equiv 1 \mod
\II$. (To see this, let $t \in \*A$, such that $t b \equiv 1 \mod \II$, and
then replace $a_1$ with $a_1 + (1 - a_1) t b$.) 
 \end{itemize}
 We have
 \begin{align*}
 \Men{bq}{a_1a_2}{} \Men{b q}{a_2}{}^{-1} 
 &= \Men{a_2 b
q (1-a_1)}{1 + a_2 d(a_1-1)}{} 
 &&\text{\pref{MenMultInv}}
 \\&= \Men{b q(1-a_1)}{1 + a_2 d(a_1-1)}{}
 &&\text{(by assumption, since $1 - a_1 \in \II$)}
 \\&= \Men{b q (1-a_1)}{a_1}{}
 &&\text{($a_2 d = 1 + bqc \equiv 1 \mod bq$)}
 \\&= \Men{b q}{a_1}{} 
 &&\text{(\fullref{BasicMennicke}{b(1-a)})}
 . \end{align*}

\begin{case} \label{3.16+}
 The general case.
 \end{case}
 Because there is no harm in interchanging $a_1$ with~$a_2$, we may assume
it is $\Men{bq}{a_2}{\Iq}$ that is equal to~$1$.

 Let $T = \mat{-1}{0}{0}{1}$. It is not difficult to see that the hypotheses
of the proposition are satisfied with $T^{-1} N T$ in the place of~$N$
(because conjugation by~$T$ is an automorphism that fixes $\Elem(2,A)$,
$\SL(2,A;\Iq)$, and $\Enorm(2,A;\Iq)$, and because $(a,-b) \in W(\II)$ for
all $(a,b) \in W(\II)$.) Therefore, the hypotheses are also satisfied with
$(T^{-1} N T) \cap N$ in the place of~$N$, so we may assume that $T$
normalizes~$N$. Then conjugation by~$T$ induces an automorphism of~$C$, so
 $$ \Men{bq}{a_2}{} = 1 = T^{-1} 1 T = T^{-1} \Men{bq}{a_2}{} T =
\Men{-bq}{a_2}{} .$$

 Adding a multiple of $a_1 a_2$ to~$b$ does not change any of the terms in
\pref{PrincMennickeSL2PF-MS2b}, so, since
$\GEN(2,1)$ holds in~$A$, we may assume that $\U{b A}/\U{b A}^2$ is cyclic.

Let $\mat{a_2}{-bq}{c}{d} \in \SL(2,A;\Iq)$. Then, by assumption and by the
formula for the inverse of a $2 \times 2$ matrix, we have
 $$\Men{b q}{a_2}{\Iq}^{-1} = \Men{- b q}{a_2}{\Iq}^{-1} = \Men{b q}{d}{\Iq} .$$

 If either $a_1$ or~$a_2$ is a square mod $b q A$, then 
 $\Men{b q}{a_1}{} \Men{b q}{a_2}{} = \Men{b q}{a_1 a_2}{}$ by Case~\ref{3.15}.

If not, then $a_1 a_2$ is a
square mod $b q A$, so, appealing to Case~\ref{3.15} again, we have
 $$\Men{b q}{a_1 a_2}{\Iq} \Men{b q}{a_2}{\Iq}^{-1} 
 = \Men{b q}{a_1 a_2}{\Iq} \Men{b q}{d}{\Iq}
 = \Men{b q}{a_1 a_2 d}{\Iq}
 = \Men{b q}{a_1}{\Iq} $$
 (because $a_2 d \equiv a_2 d + bq c = 1 \mod bq$).
 \end{proof}

\subsection{Finiteness of $\SL(2,A;\Iq)/\Enorm(2,A;\Iq)$}
\label{SL2/E2finiteSect}

We now prove the following theorem.

\begin{thm}[{(cf.\ \cite[(3.19)]{CKP})}] \label{SL2/E2Finite}
 Suppose
 \begin{itemize}
 \item $\rr$, $\xx$, $\ell$, $\ttt$, and~$\zz$ are positive integers,
  \item $A$ is an integral domain satisfying
 \begin{itemize}
 \item the stable range condition $\SRquot$,
 \item the few generators properties $\GEN(2,1)$ and $\GEN(\ttt,\rr)$,
 \item the exponent property $\EXP(\ttt,\ell)$,
 \item the unit property $\UNIT(1,\xx)$,
 and
 \item the conjugation property $\CONJ(\zz)$, 
 and
 \end{itemize}
 \item $\Iq$ is any nonzero ideal in~$A$.
 \end{itemize}
 Then $\SL(2,A;\Iq)/\Enorm(2,A;\Iq)$ is finite.
 \end{thm}

\begin{rem}
 The proof will show that the order of the quotient group is bounded by
$\ttt^\rr$.
 \end{rem}

\begin{assump}
 Throughout \S\ref{SL2/E2finiteSect}, $\rr$, $\xx$, $\ell$, $\ttt$, $\zz$,
$A$, and~$\Iq$ are as in the statement of Theorem~\ref{SL2/E2Finite}.
 \end{assump}

\begin{notation} 
  For $(a,b) \in W(\Iq)$, we set
 $$\Men{b}{a}{\Iq} =  \mat{a}{b}{*}{*} \Enorm(2,A;\Iq) \in
\SL(2,A;\Iq)/\Enorm(2,A;\Iq).$$
 \end{notation}

The key to the proof is showing that $\Men{\ }{\ }{\Iq}$ is a well-defined
Mennicke symbol. For this, we use Proposition~\ref{MennInSL2}, so it suffices to
show that
 $$\bigl[ \Elem(2,A), \SL(2,A;\Iq) \bigr] \subseteq
\Enorm(2,A;\Iq)$$
 and 
 that $\displaystyle \*{\Men{b y^2}{a}{\Iq}} = \*{\Men{b}{a}{\Iq}}$ for all
$(a,b) \in W(\II)$ and all nonzero $y \in A$. These assertions are
established in \pref{3.10} and~\pref{3.11}, respectively.

Combining Vaserstein's identity \fullref{VasId}{A4} with the conjugation
property $\CONJ(\zz)$ yields the following lemma.

\begin{lem} \label{4.7easy}
 There is a nonzero ideal~$\Iq'$ of~$A$, such that
 $$E_{1,2}^{-1} \LU(2,\Iq') \, E_{1,2} \subseteq
\gennum{\LU(2,\Iq)}{50\zz} .$$
 \end{lem}

\begin{proof}[{Proof {\rm (\cite[Lem.~4]{Vaserstein-SL2}, \cite[(4.7)]{CKP})}}]
 Fix some nonzero $q \in \Iq$. For convenience, let
 $$\conje{\LU} = E_{1,2} \LU(2,A; q A) E_{1,2}^{-1}$$
 and
 $$ \text{$\conje{E_{2,1}}(a) = E_{1,2} \, E_{2,1}(a) E_{1,2}^{-1}$, for $a \in
A$.} $$
 From the unit property \fullref{3.1}{1}, there is a unit~$u$ in~$A$, such
that $u^2 \neq 1$.
 Since $A$ satisfies $\CONJ(\zz)$ (and $E_{1,2}$ normalizes $\{E_{1,2}(*)\}$),
 $$ \text{it suffices to show $E_{2,1} \bigl( -(u^2 -1) y \bigr) \in
\gennum{\conje{\LU}}{50}$, for every $y \in \lowset_{q A}$.}$$

From \fullref{VasId}{A2} and \fullref{VasId}{A3}, we see, for any unit~$v$,
that
 $$ H(v)^{-1} \conje{E_{2,1}}(*) H(v)
 = E_{1,2}(*) \,  \conje{E_{2,1}}(*)
\, E_{1,2}(*)
 \in \gennum{ \conje{\LU} }{3} .$$
 Since $H(v)$ normalizes $\{ E_{1,2}(*)\}$, this implies that
 \begin{equation} \label{4.7easy-Hstretch}
 \text{$H(v)^{-1} \gennum{ \conje{\LU} }{j} H(v) \subseteq
\gennum{ \conje{\LU} }{3j}$ for all~$j$.}
 \end{equation}

 Because $y \in \lowset_{q A}$, we have $y \in q A$, and there exist $z \equiv
\pm 1 \mod qA$ and units $u_1$ and~$u_2$, such that $1 + yzu_1^2 = u_2^2$. By
replacing $y$ with~$-y$ if necessary, let us assume $z \equiv 1 \mod qA$.
 It is obvious that $u^2 - 1$ is a multiple
of $1 + y z u_1^2$ (since everything is a multiple of any unit). In the
notation of \fullref{VasId}{A4}, with $y' = y u_1^2$ in the role of~$y$, we
have
 $$ \text{$M(y',z) \in \gennum{ \conje{\LU}}{2}$
 and
 $ E_{1,2}(c) M(y',z)^{-1} \in \gennum{ \conje{\LU}}{3}$,} $$
 so
 $$ E_{2,1}(-w y')
 \in \gennum{ \conje{\LU}}{2 + 3 \times 3}
 = \gennum{ \conje{\LU} }{11} .$$
 Since 
 $$-w y'
 = - \frac{(u^2 - 1) y u_1^2}{1 + yz u_1^2}
 = - \frac{(u^2 - 1) y u_1^2}{u_2^2} ,$$
 conjugating by $H(u_2/u_1)$ yields the conclusion that
 \begin{align*}
 E_{2,1} \bigl( -(u^2 - 1) y \bigr)
 &= H(u_2/u_1)^{-1} E_{2,1}(-w y') \, H(u_2/u_1)
 \\&\in \gennum{ \conje{\LU} }{3 \times 11} 
 \subseteq \gennum{ \conje{\LU} }{50} ,
 \end{align*}
 as desired. 
 \end{proof}

\begin{cor}[{\cite[(3.3)]{CKP}}] \label{3.3}
 For each element~$T$ of\/ $\GL(2,K)$, there is a nonzero ideal~$\Iq'$
of~$A$, such that $T^{-1} \LU(2,\Iq') \, T \subseteq
\gennum{\LU(2,\Iq)}{(50\zz)^2}$.
 \end{cor}

\begin{proof}
 Any matrix in $\GL(2,K)$ is a product involving only diagonal matrices, the
permutation matrix $\mat{0}{1}{1}{0}$, and the elementary matrix~$E_{1,2}$,
with the elementary matrix appearing no more than twice. (This is a
consequence of the ``Bruhat decomposition.") 
 \begin{enumerate}
 \item For a diagonal matrix $T = \mat{a}{0}{0}{b}$, we have 
 $$ T^{-1} \LU(2, a^{-1} b \Iq \cap b^{-1} a \Iq \cap A ) \, T \subseteq
\LU(2,\Iq) .$$
 \item Conjugation by the permutation matrix $\mat{0}{1}{1}{0}$
interchanges $E_{1,2}(*)$ with $E_{2,1}(*)$, so $\LU(2,\Iq)$ is
invariant.
 \item For $E_{1,2}$, see \pref{4.7easy}.
 \end{enumerate}
 \end{proof}

Lemma~\ref{3.3} can be restated very cleanly in the terminology of nonstandard analysis:

\begin{cor}[{\cite[(3.6)]{CKP}}] \label{3.6}
 $\GL(2,K)$ normalizes $\Elem(2, \II)$.
 \end{cor}

\begin{proof}
 For any $T \in \GL(2,K)$, Lemma~\ref{3.3} implies that 
 $$T^{-1} \LU(2,\II) \, T \subseteq \gennum{\LU(2,\*\Iq)}{(50\zz)^2} .$$
  Since $\Iq$ is an arbitrary ideal of~$A$, we conclude, from
 polysaturation, that 
 $$ T^{-1} \LU(2,\II) \, T \subseteq \gennum{\LU(2,\II)}{(50\zz)^2}
  \subseteq \Elem(2,\II) ,$$
 as desired.
 \end{proof}

Hence, the action of $\GL(2,K)$ on $\SL(2, \*A;\II)$ induces an action on
the coset space
 $\SL(2, \*A;\II) / \Elem(2,\II)$.
 It can be shown that this is a trivial action of $\GL(2,K)$  \see{3.18},
but we now establish this only for $\SL(2,K)$. 

\begin{prop}[{\cite[(3.7)]{CKP}}] \label{3.7}
 $[\SL(2,K) ,\SL(2,\*A;\II) ] \subseteq \Elem(2,\II)$.
 \end{prop}

\begin{proof}
 Let $T$ be an arbitrary element of $\SL(2, \*A;\II)$.
 Applying \pref{3.3}, with $\*A$
in the role of~$A$, yields a nonzero ideal $\II' \subseteq \II$, such that 
 \begin{equation} \label{3.7Pf-ConjInE}
 T^{-1} \Elem(2,\II') \, T \subseteq \Elem(2,\II). 
 \end{equation}
We may assume $\II'$ is principal, by passing to a smaller ideal.
 We may write $T = X E$, with $X \in \SL(2,\*A; \II')$ and $E \in
\Elem(2,\II)$ (by Vaserstein's Lemma~1 \pref{VasersteinLemma} and the fact
that $\SL(2,\*A; \II) = \Vas(2,\*A; \II)$ \see{3.5}).
 Let $u_0$ be a unit in~$A$ satisfying the unit property \fullref{3.1}{3}
(with $\rrr = 1$), so there exist $E_1,\ldots,E_\xx \in \LU(2,\II')$,
such that
 $$ H(u_0)^{-1} X H(u_0)
 = E_1 X E_2 \cdots E_\xx
 \in E_1 T  \Elem(2,\II)
 = T \Elem(2,\II) .$$
 Then
 $$ H(u_0)^{-1} T H(u_0)
 = \bigl( H(u_0)^{-1} X H(u_0) \bigr) \bigl( H(u_0)^{-1} E H(u_0) \bigr)
 \in T \Elem(2,\II). $$
 Hence, $H(u_0)$ is in the kernel of the action on $\SL(2, \*A;\II) /
\Elem(2,\II)$. Since $\SL(2,K)$ is the smallest normal subgroup of $\GL(2,K)$
containing $H(u_0)$, this implies that all of $\SL(2,K)$ is in the kernel.
 \end{proof}

 \begin{cor}[{\cite[(3.10)]{CKP}}] \label{3.10}
 We have
 \begin{enumerate}
 \item \label{3.10-commutator}
 $\bigl[ \Elem(2,A), \SL(2,A;\Iq) \bigr] \subseteq \Enorm(2,A;\Iq)$,
 and
 \item \label{3.10-normal}
 $\Enorm(2,A;\Iq)$ is normal in $\SL(2,A;\Iq)$.
 \end{enumerate}
 \end{cor}

\begin{proof}
 \pref{3.10-commutator} We have $\SL(2,\*A;\*\Iq) = \Enorm(2,\*A;\*\Iq) \,
\SL(2,\*A;\II)$
 \see{VasLem1Normal}.
 For $T \in \Elem(2,A)$, we have
 \begin{align*}
 [T,\*\SL(2,A; \Iq) ] 
 &\subseteq [\Elem(2,A),\SL(2, \*A; \*\Iq) ] 
 \\&= [\Elem(2,A),\Enorm(2, \*A; \*\Iq) ] [\Elem(2,A), \SL(2, \*A;\II) ] 
 \\&\subseteq \Enorm(2, \*A; \*\Iq) ] \Elem(2,\II) 
 && \text{(by \ref{3.7})}
 \\&= \Enorm(2, \*A; \*\Iq) 
 \\&\subseteq \*\Enorm(2,A;\Iq)
 . \end{align*}
 By Leibniz' Principle, then $[T, \SL(2,A;\Iq)] \subseteq \Enorm(2,A;\Iq)$. This
completes the proof of the first half of the corollary.

\pref{3.10-normal} Since $\Enorm(2,A;\Iq) \subseteq \Elem(2,A)$,
part~\pref{3.10-commutator} implies 
 $$[\Enorm(2,A;\Iq), \SL(2,A;\Iq)] \subseteq \Enorm(2,A;\Iq) ,$$
 so $\Enorm(2,A;\Iq)$ is normal.
 \end{proof}

\begin{cor}[{\cite[(3.11)]{CKP}}] \label{3.11}
 If $y$ is a nonzero element of~$A$, and $(a,b)\in W(\II)$, then 
 $$ \Men{by^2}{a}{\II} = \Men{b}{a}{\II} .$$
 \end{cor}

\begin{proof}
 Because $\Elem(2, \II) \subseteq \Enorm(2, \*A;\II)$,
 we see, from \pref{3.7}, that  
 $$[\SL(2,K), \SL(2,\*A;\II)] \subseteq
\Enorm(2,\*A;\II) .$$
 Therefore
 \begin{align*}
 \Men{b y^2}{a}{\II}
 &\equiv \mat{a}{by^2}{*}{*}
 = \mat{y}{0}{0}{y^{-1}} \mat{a}{b}{*}{*} \mat{y^{-1}}{0}{0}{y} 
 \\&\equiv \mat{a}{b}{*}{*}
 \equiv \Men{b}{a}{\II}
 \mod \Enorm(2,\*A;\II)
 . \end{align*}
 \end{proof}

\begin{proof}[\bf Proof of Theorem~\ref{SL2/E2Finite}]
 We have 
 \begin{itemize}
 \item $\bigl[ \Elem(2,A), \SL(2,A;\Iq) \bigr] \subseteq
\Enorm(2,A;\Iq)$ \see{3.10}
 and 
 \item $\*{\Men{b y^2}{a}{\Iq}} = \*{\Men{b}{a}{\Iq}}$ for all $(a,b) \in
W(\II)$ and all nonzero $y \in A$ \see{3.11},
 \end{itemize}
 so Proposition~\ref{MennInSL2} implies that $\Men{\ }{\ }{\Iq}$ is a well-defined
Mennicke symbol. Therefore, its range $\SL(2,A;\Iq)/ \Enorm(2,A;\Iq)$ is
isomorphic to a quotient of the universal Mennicke group $C(\Iq)$, so the
desired conclusion is immediate from Theorem~\ref{1.8}.
 \end{proof}

\begin{rem}[{\cite[(3.18)]{CKP}}] \label{3.18}
 For $T = \mat{y^{-1}}{0}{0}{1}$, the fact that $\Men{\ }{\ }{\II}$ is a
Mennicke symbol (cf.\ proof of Theorem~\ref{SL2/E2Finite}) implies
 \begin{align*}
 T^{-1} \Men{b}{a}{\II} T
 &= \Men{b y}{a}{\II}
 = \Men{by(1-a)}{a}{\II}
 = \Men{b}{a}{\II} \Men{y(1-a)}{a}{\II}
 \\&= \Men{b}{a}{\II} \left( T^{-1} \Men{a-1}{a}{\II} T \right)
 =  \Men{b}{a}{\II} ( T^{-1} \, 1 \, T )
 = \Men{b}{a}{\II} 
 , \end{align*}
 so $T$
 acts trivially on $\SL(2,\*A;\II)/\Elem(2,\II)$.
 Combining this with \pref{3.7} yields the conclusion that all of $\GL(2,K)$
acts trivially.
 \end{rem}

\subsection{Bounded generation in $\SL(2,\BS)$}

We now deduce Theorem~\ref{EijBddGen} under the assumption that $n = 2$. (See
Theorem~\ref{BddGen>3} for the case $n \ge 3$.)

\begin{thm}[{(cf.\ \cite[(3.19)]{CKP})}] \label{BddGen2}
 Let
 \begin{itemize}
 \item $K$ be an algebraic number field,
 \item $k$ be the degree of~$K$ over~$\rational$,
 \item $B$ be an order in~$K$,
 \item $S$ be a multiplicative subset of~$B$,
 and
 \item $\Iq$ be an ideal in~$\BS$.
 \end{itemize}
 If $\BS$ has infinitely many units, then:
 \begin{enumerate}
 \item \label{BddGen2-LU}
 $\LU(2,\BS)$ boundedly generates $\Elem(2,\BS)$,
 and
 \item \label{BddGen2-LUnorm}
 the set $\LUnorm(2,\BS;\Iq)$ boundedly generates $\Enorm(2,\BS;\Iq)$.
 \end{enumerate}
 More precisely, there is a positive integer~$r$, depending only on~$n$
and~$k$, such that 
 $$ \text{$\gennum{\LU(2,\BS)}{r} = \Elem(2,\BS)$ and $\gennum{\LUnorm(2,\BS)}{r}
= \Enorm(2,\BS;\Iq)$.}$$
 \end{thm}

To establish the above result, note that Theorem~\ref{SL2/E2Finite} applies to
the above situation (by \ref{4.4}, \ref{4.5}, \ref{4.6}, and \ref{4.7}), so the
desired conclusion follows from the Compactness Theorem \see{CpctnessCor}.

\section{Bounded generation of normal subgroups} \label{NormalSect}

\begin{thm}[{(cf.\ \cite[(2.7) and (3.21)]{CKP})}] \label{BddGenNormal}
 Let
 \begin{itemize}
 \item $n$ be a positive integer,
 \item $K$ be an algebraic number field,
 \item $k$ be the degree of~$K$ over~$\rational$,
 \item $B$ be an order in~$K$, and
 \item $S$ be a multiplicative subset of~$B$.
 \end{itemize}
 Assume that either $n \ge 3$ or $\BS$ has infinitely many units.
 \begin{enumerate}
 \item  \label{BddGenNormal-normal}
 If $\genset^{\normal}$ is any subset of\/ $\SL(n,\BS)$, such that $g^{-1}
\genset^{\normal} g = \genset^{\normal}$, for every $g \in \Elem(n,\BS)$
{\rm(}and $\genset^{\normal}$ does not consist entirely of scalar
matrices{\rm)}, then $\genset^{\normal}$ boundedly generates a finite-index
subgroup of\/ $\SL(n,\BS)$. 
 \item  \label{BddGenNormal-Gamma}
 For any finite-index subgroup~$\Gamma$ of\/ $\SL(n,\BS)$, the
set\/ $\LU(n, \BS) \cap \Gamma$ of elementary matrices in\/~$\Gamma$ boundedly
generates a subgroup of finite index in\/~$\Gamma$.
 \end{enumerate}
 \end{thm}

\begin{rem} \label{normalbound} 
 In the situation of part~\pref{BddGenNormal-normal} of the above theorem, we
have $\gennum{\genset^{\normal}}{r} = \gen{\genset^{\normal}}$, for some~$r$
that depends on~$k$, $n$, $\#(A/\Iq)$, and the minimal number of generators
of~$\Iq$, where $\Iq = \Iq(\genset^{\normal})$~is a certain ideal defined in
the statement of Proposition~\ref{Normal} below. (The minimal number of generators
of~$\Iq$ is certainly finite, since $\BS$ is Noetherian. In the situation of
Theorem~\ref{ConjBddGen} of the introduction, the minimal number of generators of
$\Iq(\genset^{\normal})$ is bounded by $n^2 \cdot \#\genset$.)
 \end{rem}

\begin{rem}
 We will use the Compactness Theorem to establish Theorem~\ref{BddGenNormal}, but
a more straightforward proof can be obtained by applying nonstandard analysis.
All of the cases are very similar, so let us describe only the proof of
\fullref{BddGenNormal}{normal} when $n \ge 3$. There is some nonzero
(principal) ideal~$\Iq'$ of~$\*\BS$, such that $\gen{\*\genset^{\normal}}$
contains $\Enorm(n,A;\Iq')$ (see Theorem~\ref{Sandwich} below).
 From Theorem~\ref{SLnq/EnqFinite}, we know that 
 $\SL(n,\*A; \Iq')/\Enorm(n,\*A;\Iq')$  is finite.
 Furthermore, since $\SL(n,A)/\SL(n,A;\Iq)$ is finite, for every nonzero
ideal~$\Iq$ of~$A$, Leibniz' Principle implies that
 $\SL(n,\*A)/\SL(n,\*A; \Iq')$ is $*$-finite.
 Therefore, the coset space $\SL(n,\*A)/\Enorm(n,\*A;\Iq')$ is $*$-finite,
which means there is a $*$-finite set~$\Omega$, such that
$\Enorm(n,\*A;\Iq') \Omega = \SL(n,\*A)$. Then
 $$ \*{\gen{\genset^{\normal}}}
 \subseteq \SL(n,\*A)
 = \Enorm(n,\*A;\Iq') \Omega
 \subseteq \gen{\*\genset^{\normal}} \Omega ,$$
 so the desired bounded generation follows from ($\ref{2.1-4} \Rightarrow
\ref{2.1-1}$) of Proposition~\ref{2.1}.
 \end{rem}

\subsection{Part \ref{BddGenNormal-normal} of Theorem~\ref{BddGenNormal}}

When $n \ge 3$, the proof of \fullref{BddGenNormal}{normal} is based on the
following description of normal subgroups of $\SL(n,A)$.

\begin{thm}[{(Sandwich Condition
\cite[Thm.~4.2e]{Bass-KTheoryStable}, \cite[4.2.9, p.~155]{HahnOMeara})}] 
\label{Sandwich}
 Suppose
 \begin{itemize}
 \item $A$ is a commutative ring that satisfies the stable range condition
$\SR_2$,
 \item $n \ge 3$,
 and
 \item $N$ is a subgroup of\/ $\SL(n,A)$ that is normalized by $\Elem(n,A)$.
 \end{itemize}
 Then there is an ideal~$\Iq$ of~$A$, such that
 \begin{enumerate}
 \item $N$ contains $\Enorm(n,A;\Iq)$,
 and
 \item each element of~$N$ is congruent to a scalar matrix, modulo~$\Iq$.
 \end{enumerate}
 \end{thm}

As a replacement for the Sandwich Condition when $n = 2$, we have the
following elementary observation, essentially due to Serre.

\begin{lem}[{(\cite[Lem.~1.3]{CostaKeller}, cf.
\cite[pp.~492--493]{Serre-CSPSL2})}] \label{u4inN}
 Suppose
 \begin{itemize}
 \item $A$ is a commutative ring,
 \item $N$ is a subgroup of\/ $\SL(2,A)$ that is normalized by $\Elem(2,A)$,
 \item $\mat{a}{b}{c}{d}$ is any element of~$N$,
 and
 \item $u$ is a unit in~$A$, such that $u^2 \equiv 1 \mod cA$.
 \end{itemize}
 Then $\Enorm \bigl( 2,A; (u^4-1) A \bigr) \subseteq N$.
 \end{lem}

Now \fullref{BddGenNormal}{normal} is obtained by applying the Compactness
Theorem \see{CpctnessCor} to the following proposition. 

\begin{defn}
 The \emph{level ideal} of a subset~$\genset$ of $\SL(n,A)$ is the smallest
ideal~$\Iq$ of $\SL(n,A)$, such that the image of $\genset$ in $\SL(n,A;\Iq)$
consists entirely of scalar matrices.
 \end{defn}

\begin{prop} \label{Normal}
 Suppose
 \begin{itemize}
 \item $\jj,\mm,\rr,\ell,\rrr,\zz$ are positive integers,
 \item $A$ is an integral domain satisfying $\SRquot$, $\GEN(2,1)$,
$\GEN(2\ttt,\rr)$, and $\EXP(2\ttt,\ell)$,
 \item $\genset^{\normal}$ is any {\rm(}nonempty{\rm)} subset of\/ $\SL(n,A)$, such
that $g^{-1} \genset^{\normal} g = \genset^{\normal}$, for every $g \in \Elem(n,A)$ {\rm(}and
$X$ does not consist entirely of scalar matrices{\rm)},
 \item either $n \ge 3$ and
 \begin{itemize}
 \item $\Iq$ is the level ideal of $\genset^{\normal}$,
 and
 \item there is a $\jj$-element subset $\genset_0$ of~$\genset^{\normal}$,
such that $\Iq$ is the level ideal of~$\genset_0$,
 \end{itemize}
 \item or $n = 2$ and
 \begin{itemize}
 \item $A$ satisfies $\UNIT(\rrr,\xx)$ and $\CONJ(\zz)$,
 \item $\mat{a}{b}{c}{d} \in \genset^{\normal}$,
 \item $u$ is a unit in~$A$, such that $u^2 \equiv 1 \mod cA$,
 and
 \item $\Iq = (u^4 - 1) A$
 \end{itemize}
 and
 \item $\#(A/\Iq) \le \mm$.
 \end{itemize}
 Then $\genset^{\normal}$ generates a finite-index subgroup of\/ $\SL(n,A)$. 
 \end{prop}

\begin{proof}
 Since $\gen{\genset^{\normal}}$ is obviously normalized by $\Elem(n,A)$, we know, from
\pref{Sandwich} or \pref{u4inN}, that $\gen{\genset^{\normal}}$ contains
$\Enorm(n,A;\Iq)$. Thus, 
 $$ \# \frac{\SL(n,A)}{\gen{\genset^{\normal}}} 
 \le \# \frac{\SL(n,A)}{\Enorm(n,A;\Iq)}
 \le \# {\SL(n,A/\Iq)}
 \cdot \# \frac{\SL(n,A;\Iq)}{\Enorm(n,A;\Iq)} 
 < \infty ,$$
 by \pref{SLnq/EnqFinite} or \pref{SL2/E2Finite}.
 \end{proof}

\subsection{Part~\ref{BddGenNormal-Gamma} of Theorem~\ref{BddGenNormal}}

In preparation for the proof of \fullref{BddGenNormal}{Gamma}, we establish
some preliminary results. (We need only the corollary that follows.)
 The following theorem is only a special case of a result that is valid for
all Chevalley groups, not only $\SL(n,A)$.

\begin{thm}[{(Tits \cite[Prop.~2]{Tits-SysGenGrpCong})}] \label{Tits-EnormInE}
 If $n \ge 3$ and $\Iq$ is any ideal of any commutative ring~$A$, then
$\Enorm(n,A;\Iq^2) \subseteq \Elem(n,\Iq)$.
 \end{thm}

See Definition~\ref{VasDefn} for the definition of $\Vas(2,A;\Iq)$.

\begin{lem}[{\cite[(3.8) and (3.20)]{CKP}}] \label{SL2Lemma}
 Suppose
 \begin{itemize}
 \item $A$ is as in the statement of Theorem~\ref{SL2/E2Finite},
 and
 \item $\Iq$ is a nonzero ideal in~$A$.
 \end{itemize}
 Then:
 \begin{enumerate}
 \item \label{SL2Lemma-3.8}
 $\Elem(2,\Iq)$ is a normal subgroup of\/ $\Vas(2,A;\Iq)$. 
 \item \label{SL2Lemma-3.20}
 If $A/\Iq^2$ is finite, then $\Elem(2,\Iq)$ contains $\Enorm(2,A; \Iq')$, for
some nonzero ideal~$\Iq'$.
 \end{enumerate}
 \end{lem}

\begin{proof}
 We use nonstandard analysis. (See Notation.~\ref{ODefn} for the definition of
the ideal~$\II$.)

 \pref{SL2Lemma-3.8}
 Let $T \in \Vas(2,A;\Iq)$ and $E_1 \in \LU(2,\Iq)$. By Vaserstein's Lemma~1
\pref{VasersteinLemma}, we may write $T = X E$ with $X \in \Vas(2,\*A;\II)$
and $E \in \Elem(2,\*\Iq)$.
 Then 
 $$[T,E_1] = [XE,E_1] = E^{-1}[X,E_1]E_1^{-1}EE_1 .$$
 It is obvious that $E_1,E \in \Elem(2,\*\Iq)$, and, by \pref{3.7}, we
have 
 $$[X,E_1] \in \Elem(2,\II) \subseteq \Elem(2,\*\Iq) .$$
 Hence 
 $$[T,E_1] \in \Elem(2,\*\Iq) \subseteq \*\Elem(2,\Iq) .$$
 By Leibniz' Principle, $[T,E_1] \in \Elem(2,\Iq)$.

\pref{SL2Lemma-3.20} 
 Let $w_1,\ldots,w_r$ be coset representatives for $\SL(2,A;\Iq^2)$ in
$\SL(2,A)$.
 For each~$i$, we have 
 $\LU(2,q_i A) \subseteq w_i^{-1} \Elem(2,\Iq)w_i$, for some nonzero $q_i$
\see{3.3}, so 
 $$ \text{$\displaystyle H = \bigcap_{i=1}^r \left( w_i^{-1} \Elem(2,\Iq)w_i
\right)$ contains $\LU(2, q_1 q_2\cdots q_r A)$.} $$
 From \pref{SL2Lemma-3.8} (and because $\SL(2,A;\Iq^2) \subseteq
\Vas(2,A;\Iq)$), we see that $H$ is the intersection of all of the conjugates
of $\Elem(2,\Iq)$, so $H$ is normal. Therefore $H$ contains $\Enorm(2,A; q_1
q_2\cdots q_r A)$.
 \end{proof}

\begin{cor} \label{BddGenContEnorm}
 Assume the situation of Theorem~\ref{BddGenNormal}.
 If\/ $\Iq$ is any nonzero ideal of $\BS$, then there exist a nonzero
ideal\/~$\Iq'$ of $\BS$ and a positive integer~$r$, such that $\Enorm(n,\BS;
\Iq') \subseteq \gennum{\LU(n,\Iq)}{r}$.
 \end{cor}

\begin{proof}[Proof {\rm (sketch)}]
 We apply a compactness argument to Lemma~\fullref{SL2Lemma}{3.20} (if $n = 2$)
or Theorem~\ref{Tits-EnormInE} (if $n \ge 3$). These results show (under
appropriate hypotheses) that there is an ideal $\Iq'$ of~$A$, such that
 $$\Enorm(n,A;\Iq') \subseteq \Elem(n,\Iq) = \gen{\LU(n,\Iq)} .$$
 By bounded generation of $\Enorm(n,A; \Iq')$ and $\Elem(n,A)$
\cfand{BddGen>3}{BddGen2}, there is some positive integer $r_0$, such that
 $$ \Enorm(n,A;\Iq') = \gennum{
 \bigcup_{\textstyle E \in \gennum{\LU(n,A)}{r_0}} E^{-1} \LU(n, \Iq') E
 }{\raise 10pt \hbox{$r_0$}}
 .$$
 Thus, the desired result is a consequence of the Compactness Theorem
\pref{Cpctness}.
 \end{proof}

\begin{proof}[\bf Proof of Theorem~\fullref{BddGenNormal}{Gamma}]
 Because $\Gamma$ has finite index, there is some nonzero ideal~$\Iq$ of
$\BS$, such that $\Enorm(n, \BS; \Iq) \subset \Gamma$.
 From Theorem~\ref{SL2/E2Finite} (if $n = 2$) or Theorem~\ref{SLnq/EnqFinite} (if $n
\ge 3$), and the fact that $\SL(n,\BS)/\SL(n,\BS;\Iq)$ is finite, we see that
 \begin{equation} \label{SLn/EnormFinite}
 \text{$\Enorm(n,\BS;\Iq)$ is a subgroup of finite index in $\SL(n,\BS)$.}
 \end{equation}
 Thus, we may assume $\Gamma = \Enorm(n, \BS; \Iq)$, so 
 $$\LU(n, \BS) \cap \Gamma = \LU(n, \Iq) .$$
 We have $\Enorm(n,\BS;\Iq') \subseteq \gennum{\LU(n,\Iq)}{r}$, for
some nonzero ideal $\Iq'$ of~$\BS$ and some positive integer~$r$
\see{BddGenContEnorm}. Since $\Enorm(n,\BS;\Iq')$ has finite index in
$\SL(n,\BS)$ \see{SLn/EnormFinite}, this implies $\gen{\LU(n,\Iq)} =
\gennum{\LU(n,\Iq)}{r'}$, for some positive integer~$r'$.
 \end{proof}

For $n \ge 3$, the bounded generation of normal subgroups also remains
valid when the group $\SL(n,\BS)$ is replaced by a subgroup of finite index
\see{BddGenGammaN}.

\begin{thm}[{(Bak \cite[Cor.~1.2]{BakNormalized})}]
 Let
 \begin{itemize}
 \item $n \ge 3$,
 \item $A$ be a commutative ring satisfying the stable range condition
$\SR_2$,
 \item $\Iq$ be a nonzero ideal of~$A$,
 and
 \item $N$ be a noncentral subgroup of\/ $\SL(n,A)$.
 \end{itemize}
 If $N$ is normalized by $\Enorm(n,A;\Iq)$, then $N$ contains
$\Enorm(n,A;\Iq')$, for some nonzero ideal~$\Iq'$ of~$A$.
 \end{thm}

\begin{cor} \label{BddGenGammaN}
 Let
 \begin{itemize}
 \item $n \ge 3$, 
 \item $K$, $k$, $B$, and~$S$ be as in Theorem~\ref{BddGenNormal},
 \item $\Gamma$ be any subgroup of finite index in $\SL(n,\BS)$,
 and
 \item $\genset^{\normal}$ be is any subset of\/ $\Gamma$, such that $g^{-1}
\genset^{\normal} g = \genset^{\normal}$, for every $g \in \Gamma$ {\rm(}and
$\genset^{\normal}$ does not consist entirely of scalar matrices{\rm)}.
 \end{itemize}
 Then $\genset^{\normal}$ boundedly generates a finite-index subgroup
of\/~$\Gamma$. 
 \end{cor}

\begin{proof}
 Applying a compactness argument (as in the proof of
Corollary~\ref{BddGenContEnorm}) to the theorem yields the conclusion that there
exist a nonzero ideal\/~$\Iq'$ of $\BS$ and a positive integer~$r$, such that
$\Enorm(n,\BS; \Iq') \subseteq \gennum{\genset^{\normal}}{r}$. The proof is
completed by arguing as in the final paragraph of the proof of
Theorem~\fullref{BddGenNormal}{Gamma}, with $\genset^{\normal}$ in the place of
$\LU(n,\Iq)$.
 \end{proof}


\begin{thebibliography}{}

\bibitem[AM]{AdianMennicke}
 S.~I.~Adian and J.~Mennicke,
 {\em On bounded generation of\/ $\SL_n(\integer)$},
 Internat. J. Algebra Comput. \textbf{2}  (1992),  no. 4, 357--365.

\bibitem[Bak]{BakNormalized}
 A.~Bak,
 {\em Subgroups of the general linear group normalized by relative elementary
matrices},
 in R.~K.~Dennis, ed.,
 \emph{Algebraic $K$-Theory {\rm(}Oberwolfach, 1980\/{\rm)},} Part~II,
 Springer Lecture Notes \#967.
 Springer, New York, 1982,
 pp.~1--22.
 
 \bibitem[Bar]{Bardakov}
 V.~G.~Bardakov,
 \emph{On the decomposition of automorphisms of free modules into simple
 factors},
 Izv. Math. \textbf {59} (1995) 333--351.

\bibitem[Ba1]{Bass-KTheoryStable}
 H.~Bass,
 {\em $K$-theory and stable algebra},
  Inst. Hautes \'Etudes Sci. Publ. Math. \textbf{22} (1964), 5--60.


\bibitem[Ba2]{Bass-AlgKTheoryBook}
  H.~Bass,
 \emph{Algebraic K-theory},
 Benjamin, New York, 1968.

\bibitem[BMS]{BMS}
  H.~Bass, J.~Milnor, and J.-P.~Serre,
 {\em Solution of the Congruence Subgroup Problem for $\SL_n$ {\rm(}$n \ge
3${\rm)} and $\mathop{\rm Sp}_{2n}$ {\rm(}$n \ge 2${\rm)}},
 Inst. Hautes \'Etudes Sci. Publ. Math. \textbf{33} (1967), 59--137.


\bibitem[CK1]{CarterKeller-BddElemGen}
  D.~Carter and G.~Keller,
 {\em Bounded elementary generation of\/
 $\SL_n(\mathscr{O})$},
 Amer. J. Math. \textbf{105} (1983), 673--687.

\bibitem[CK2]{CarterKeller-ElemExp}
  D.~Carter and G.~Keller,
 {\em Elementary expressions for unimodular matrices}, 
 Comm. Algebra \textbf{12} (1984), 379--389.

\bibitem[CKP]{CKP}
 D.~Carter, G.~Keller, and E.~Paige,
 {\em Bounded expressions in $\SL(n,A)$} (unpublished). 

\bibitem[CW]{CookeWeinberger}
 G.~Cooke and P.~Weinberger,
 {\em On the construction of division chains in algebraic number rings, with
applications to $\SL_2$}, 
 Comm. Algebra \textbf{3(6)} (1975), 481--524. 

\bibitem[CoK]{CostaKeller}
  D.~L.~Costa and G.~Keller,
 {\em On the normal subgroups of $\SL(2,A)$}, 
  J. Pure Appl. Algebra \textbf{53} (1988), 201--226.

\bibitem[DV]{DennisVaserstein}
 R.~K.~Dennis and L.~N.~Vaserstein,
 {\em On a question of M.~Newman on the number of commutators},
 J.~Algebra \textbf{118} (1988), 150--161.

\bibitem[ER1]{ErovenkoRapinchukCR}
 V.~Erovenko and A.~Rapinchuk,
 {\em Bounded generation of some $S$-arithmetic orthogonal groups},
 C. R. Acad. Sci. Paris S\'er.~I Math. \textbf{333}  (2001),  no. 5,
395--398.

\bibitem[ER2]{ErovenkoRapinchukPreprint}
 V.~Erovenko and A.~Rapinchuk,
 {\em Bounded generation of S-arithmetic subgroups of
isotropic orthogonal groups over number fields},
 preprint, 2005.


\bibitem[GS]{GrunewaldSchwermer}
 F.~Grunewald and J.~Schwermer,
 {\em Free non-abelian quotients of $\SL_2$ over orders of imaginary quadratic
numberfields},
 J. Algebra \textbf{69} (1981), 298--304.

\bibitem[HOM]{HahnOMeara}
 A.~J.~Hahn and O.~T.~O'Meara,
 \emph{The Classical Groups and K-Theory,}
 Springer, New York, 1989.


\bibitem[Li]{Liehl-Beschrankte}
 B.~Liehl,
 {\em Beschr\"ankte Wortl\"ange in $\SL_2$}, 
 Math. Z. \textbf{186} (1984), 509--524.

\bibitem[LM]{LoukanidisMurty}
 D.~Loukanidis and V.~K.~Murty,
 {\em Bounded generation for $\SL_n$ {\rm(}$n \ge 2${\rm)} and $\mathop{\rm
Sp}_n$
 {\rm(}$n \ge 1${\rm)}}
 (preprint).


\bibitem[Mu]{Murty-BddGen}
 V.~K.~Murty,
 {\em Bounded and finite generation of arithmetic groups},
 in: K.~Dilcher, ed.,
 \emph{Number theory {\rm(}Halifax, NS, 1994\/{\rm)},}
 CMS Conf. Proc. \#~15, Amer. Math. Soc., Providence, RI, 1995,
 pp.~249--261.

 \bibitem[Os]{Ostmann}
 H.~Ostmann,
 \emph{Additive Zahlentheorie,} Vol. 2,
 Springer-Verlag, Berlin, 1956.

\bibitem[Ra]{Rapinchuk-CSPFinWidth}
 A.~S.~Rapinchuk,
 {\em The congruence subgroup problem for arithmetic groups of finite
width} (Russian),
 Dokl. Akad. Nauk SSSR \textbf{314} (1990), 1327--1331;
 translation in 
 Soviet Math. Dokl. \textbf{42} (1991),  no. 2, 664--668.


\bibitem[Se]{Serre-CSPSL2}
 J.-P.~Serre,
 {\em Le probl\`eme des groupes de congruence pour $\SL_2$},
 Ann. Math. \textbf{92} (1970), 489--527.

\bibitem[Sh]{Shalom-BddGen}
 Y.~Shalom,
 {\em Bounded generation and Kazhdan's property~$(T)$},  
 Inst. Hautes \'Etudes Sci. Publ. Math. \textbf{90} (1999), 145--168.
 
 \bibitem[SS]{SivatskiStepanov}
 A.~S.~Sivatski and A.~V.~Stepanov,
 {\em On the word length of commutators in $\GL_n(R)$},
 K-Theory \textbf{17} (1999), 295--302.

\bibitem[SV]{StepanovVavilov}
 A.~Stepanov and N.~Vavilov,
 {\em Decomposition of transvections: a theme with variations},
 $K$-Theory \textbf{19} (2000), 109--153.

\bibitem[SL]{StroyanLuxemburg}
 K.~D.~Stroyan and W.~A.~J.~Luxemburg,
 \emph{Introduction to the Theory of Infinitesimals},
 Academic Press, New York, 1976.

\bibitem[Ta1]{Tavgen-BddGenChev}
 O.~I.~Tavgen,
 {\em Bounded generation of Chevalley groups over rings of algebraic
$S$-integers} (Russian),
 Izv. Akad. Nauk SSSR Ser. Mat. \textbf{54}  (1990),  no. 1, 97--122,
221--222; 
 translation in Math. USSR-Izv. \textbf{36} (1991),  no. 1, 101--128.

\bibitem[Ta2]{Tavgen-FinWidthArith}
 O.~I.~Tavgen,
 {\em Finite width of arithmetic subgroups of Chevalley groups of rank
$\ge 2$} (Russian),
 Dokl. Akad. Nauk SSSR \textbf{310} (1990), no. 4, 802--806;
 translation in Soviet Math. Dokl. \textbf{41} (1990), no. 1, 136--140.

\bibitem[Ti]{Tits-SysGenGrpCong}
 J.~Tits,
 {\em Syst\`emes g\'en\'erateurs de groupes de congruence},
 C.~R.~Acad. Sci. Paris, S\'er.~A  \textbf{283} (18 October 1976),
693--695.

\bibitem[vdK]{vanderKallen-SL3(C[X])}
 W.~van der Kallen,
 {\em $\SL_3\bigl(\complex[X]\bigr)$ does not have bounded word length},
 in R.~K.~Dennis, ed.,
 \emph{Algebraic $K$-Theory {\rm(}Oberwolfach, 1980\/{\rm)},} Part~I,
 Springer Lecture Notes \#966.
 Springer, New York, 1982,
 pp.~357--361.

\bibitem[Va]{Vaserstein-SL2}
 L.~N.~Vaserstein,
 {\em On the group $\SL_2$ over Dedekind rings of arithmetic type}
(Russian), Mat. Sb. (N.S.) \textbf{89(131)} (1972), 313--322, 351;
 translation in
 Math. USSR Sb. \textbf{18} (1972), 321--332.

\bibitem[Za]{Zakiryanov}
 K.~Kh.~Zakir'yanov,
 {\em Symplectic groups over rings of algebraic integers have finite width
 over the elementary matrices},
Algebra and Logic \textbf{24} (1985) 436--440.

\end{thebibliography}
\end{document}